\documentclass[12pt,a4paper]{statsoc}
\usepackage[utf8]{inputenc}
\usepackage[english]{babel}
\usepackage{amsmath}
\usepackage{amsfonts}
\usepackage{amssymb}
\usepackage{amsbsy}
\usepackage{bm}
\usepackage{graphicx}
\graphicspath{{figures/}}
\usepackage[a4paper]{geometry}
\usepackage[hidelinks]{hyperref} 
\hypersetup{
    colorlinks = false
}
\usepackage{natbib,enumerate}
\usepackage{subcaption}
\usepackage{newfloat}
\DeclareFloatingEnvironment[
    name=Table,
]{alternativetable}

\title{Stochastic differential equation based on a multimodal potential to model movement data in ecology}
\author[Gloaguen { et al.}]{Pierre Gloaguen}
\address{AgroParisTech, Paris, France.}
\address{Institut français de recherche pour l'exploitation de la mer, Nantes, France.}
\email{pierre.gloaguen@gmail.com}
\author[Gloaguen { et al.}]{Marie-Pierre Etienne}
\address{AgroParisTech, Paris, France.}
\author[Gloaguen {et al.}]{Sylvain Le Corff}
\address{Laboratoire de Math\'ematiques d'Orsay, Univ. Paris-Sud, CNRS, Universit\'e
Paris-Saclay, 91405 Orsay, France.}

\newcommand{\x}{\mathsf{x}}

\newcommand{\R}{\mathbb{R}}
\newcommand{\rmd}{\text{d}}
\newcommand{\bo}[1]{\mathbf{#1}}
\newcommand{\eqsp}{\,}
\newcommand{\eqdef}{:=}
\newcommand{\inv}{^{-1}}

\newcommand{\bs}[1]{\boldsymbol{#1}}
\begin{document}
\maketitle
\begin{abstract}
This paper proposes a new model for individuals movement in ecology. The movement process is defined as a solution to a stochastic differential equation whose drift is the gradient of a multimodal potential surface. This offers a new flexible approach among the popular potential based movement models in ecology. To perform parameter inference, the widely used Euler method is compared with two other pseudo-likelihood procedures and  with  a Monte Carlo Expectation Maximization approach based on exact simulation of diffusions. Performances of all methods are assessed with simulated data and with a data set of fishing vessels trajectories. We show that the usual Euler method performs worse than the other procedures for all sampling schemes.
\end{abstract}
 		
\section{Introduction}
Statistical inference of movement models provides many insights on the ecological features that explain population-level dynamics. 
These analyses are crucial to wildlife managers to understand complex animal behaviors \citep{Chavez+2006}. 
In fisheries science, understanding the underlying patterns responsible for spatial use of the ocean is a key aspect of a sustainable management \citep{Chang2011}. 
Both fields promote now large programs to deploy Global Positioning System (GPS) device. For instance, we might mention, among others, the Tagging of Pelagic Predators program (TOPP, \url{www.gtopp.org}), the TORSOOI (\url{www.torsooi.com}) program  for marine animals, the Elephant without Borders (\url{www.elephantswithoutborders.org/tracking.php}) or WOLF GPS (\url{www.wolfgps.com})   programs for terrestrial wildlife.
In the European Union,  since the 1st of January 2012, fishing vessels above 12m  are mandatory  equipped with a Vessel Monitoring System which has become a standard tool of fisheries monitoring worldwide. 
 As a result, such programs produce huge datasets which have been largely used to understand, explain and predict animals or vessels movements. The GPS-type loggers can be set with several acquisition frequencies which should be adjusted to the experimental setting. 
 
In the last few years, growing attention has been dedicated to continuous time and continuous space Markovian models as a realistic and flexible framework to model such data \citep{blackwell:1997,brillinger:et:al:2002,preisler:et:al:2004}. 
These papers focus on stochastic differential equations (SDE) to describe and analyze animal trajectories. Since the model based on a pure diffusion process  introduced in \cite{skellam_random_1951} a large variety of SDE has been proposed to model wildlife behavior, see \cite{preisler:et:al:2004} for references about elephant-seal or birds migrations.
In \cite{brillinger:et:al:2001a}, the authors introduces SDE based on a potential surface to capture the directional bias of animal patterns. 
This potential function is assumed to reflect the attractiveness of the environment and regions where species are likely to travel to.  
In the framework proposed by \cite{brillinger:2010}, the drift of the SDE is given by  the gradient of a potential map $P_\eta$ which depends on an unknown parameter $\eta$ while the diffusion coefficient is a smooth function. 
This flexible framework has been widely used in movement ecology for the last 20 years \citep{blackwell:1997, brillinger:et:al:2001a, brillinger:et:al:2002, brillinger:et:al:2011,harris:blackwell:2013,preisler:et:al:2004,preisler:et:al:2013}.  
In \cite{blackwell:1997}, \cite{blackwell:et:al:2015}, and \cite{harris:blackwell:2013}, the authors introduce a quadratic potential function which corresponds to a bivariate Ornstein-Uhlenbeck process to model the position of an animal. 
The Ornstein-Ulhenbeck process offers a convenient framework to represent attractiveness of a given area but it remains restrictive as animals are not prone to revert to a single attractive zone. 
This process owes its popularity to the simple form of the maximum likelihood estimator of $\eta$ and of the diffusion coefficient as its transition density is Gaussian. 
In \cite{preisler:et:al:2013}, the authors propose different potential functions, with more complex features such as multiple attractive regions, as a baseline of a flexible framework to describe animal movement. 

In this paper, we propose a new model where the potential function is a mixture of Gaussian shaped functions. This model was suggested in \cite{preisler:et:al:2013}, but, to the best of our knowledge, it has never been used in movement ecology or applied to practical cases. 
Each attractive region is modeled as a Gaussian shaped function with its own mean location and information matrix characterizing its dispersion. The drift of the SDE ruling animal movement is a mixture of these functions, with unknown weights, centers and information matrices that should be inferred from movement data. 

Designing efficient procedures to infer this model is a very challenging task as some elementary properties, such as the transition density of the solution to the SDE, are not available in closed form.
Different statistical methods have been proposed to perform maximum likelihood estimation with unknown transition densities (see \citealp{iacus2009simulation,kessler:lindner:sorensen:2012}).  
The most convenient procedures to estimate discretely observed diffusion processes rely on discretization schemes to approximate the SDE between two observations and use a surrogate of the likelihood function to compute an approximate maximum likelihood estimator of $\eta$. 
The simplest likelihood approximation is the Euler-Maruyama method which uses constant drift and diffusion coefficients  on each interval.  To the best of our knowledge, maximum likelihood inference procedures in movement ecology  in the case of unknown transition densities are only based on the Euler-Maruyama method (see \citealp{brillinger:2010, preisler:et:al:2013}). 
The maximum likelihood estimator computed with this method may be proved to be consistent under some assumptions on the total number of observations $n$ and on the time step $\Delta_n$ between observations (which is assumed to vanish as $n$ goes to infinity). 
However, there is no guarantee that using GPS tagging in ecology provides observations at a sufficiently high sampling rate and it might produce a highly biased estimator. Indeed, the chosen frequency is a crucial parameter to monitor animal behaviors over a given time horizon since GPS batteries depend on the number of recorded relocations. 
Therefore, depending on the needs of the given study, the recording frequency can be small. 
Moreover, real life experiments in complex environments (water for marine mammals for instance) might result in strongly irregular time step acquisitions. 
Therefore, alternatives to the Euler-Maruyama method should be considered to estimate parameters in potential based movement models. For instance, the bias of the Euler-Maruyama method may be reduced using higher order discretization schemes.

A particularly interesting higher order scheme is given by the Ozaki discretization which proposes a linear approximation of the drift coefficient combined with a constant approximation of the diffusion term \citep{ozaki:1992,shoji:1998}. 
More recently, \cite{kessler:1997}, \cite{kessler:lindner:sorensen:2012} and \cite{uchida:yoshida:2012} introduced another Gaussian based approximation of the transition density between two consecutive observations using higher order expansions of the conditional mean and variance of an observation given the observation at the previous time step. 
Another approach to approximate the likelihood of the observations based on Hermite polynomials expansion was introduced by \cite{ait-sahalia:1999,ait-sahalia:2002,ait-sahalia:2008} and extended in several directions recently, see \cite{li:2013} and all the references therein. However these methods also induce, theoretically, a bias due to the likelihood approximation.

To avoid this systematic bias associated with every discretization procedure, \cite{beskos:roberts:2005} and \cite{beskos:papaspiliopoulos:roberts:2006} proposed an exact algorithm to draw samples from finite-dimensional distributions of potential based diffusion processes. 
This sampling method is based on a rejection sampling step and  may be used to sample from any finite dimensional distribution of the target SDE. 
This algorithm is a pivotal tool to obtain an unbiased estimate of the intermediate quantity of the Expectation Maximization (EM) \citep{dempster:maximum:1977}. 
As a byproduct, this also allows to sample trajectories exactly distributed according to the estimated model which is not possible with the three other methods. 
However, the computational complexity of the acceptance rejection step of these exact algorithms, which requires to sample skeletons of trajectories between each pair of observations, grows with the time step between observations.

 There is no theoretical nor practical studies to choose the method which achieves the  best trade off between fast convergence and accuracy in movement ecology. We investigate the performance of four different estimation methods for the new model proposed in this paper:
\begin{enumerate}[-]
\item the Euler-Maruyama method ;
\item the Ozaki local linearization method, proposed in \cite{ozaki:1992} and \cite{shoji:1998} ;
\item the Kessler high order Gaussian approximation method, proposed in \cite{kessler:1997} and refined using an adaptive procedure in \cite{uchida:yoshida:2012} ;
\item the Exact Algorithm based Monte Carlo EM method, proposed in \cite{beskos:papaspiliopoulos:roberts:fearnhead:2006}.
\end{enumerate}
The paper is organized as follows. The potential based model is introduced in Section~\ref{sec:model} and the inference methods used in the paper are presented in Section~\ref{sec:inference}. A simulation study to evaluate the robustness of  the statistical methods with the specific potential function proposed in the paper is given in Section \ref{sec:sim:study}. In order to evaluate the impact of the sampling rate, three sampling scenarios are considered, high frequency, intermediate frequency, and low frequency. Finally, an  application of our new model to real data, together with an evaluation of the performance of each inference method, is done in Section~\ref{subsec:data:application} for the estimation of attractive zones of two different French fishing vessels.
\section{Potential based movement model}
\label{sec:model}
In this section, we propose a new model to describe the position process of an individual. The latter modelled as is a stochastic process $(\bo{X}_t)_{0\le t\leq T}$ observed between times 0 and $T$. 
$(\bo{X}_t)_{0\le t\leq T}$ is assumed  to be the solution to the following time homogeneous SDE:
\begin{equation}
\label{eq:target:sde1}
\bo{X}_0 = \bo{x}_{0}\quad\mbox{and}\quad\rmd \bo{X}_t=b_{\eta}(\bo{X}_t)\rmd t + \gamma \rmd \bo{W}_t\eqsp,
\end{equation}
where $(\bo{W}_t)_{0\le t \le T}$ is a standard Brownian motion on $\mathbb{R}^2$, and  $\gamma\in\mathbb{R}_+^\star$ is an unknown diffusion parameter. Moreover, we assume that the drift follows the gradient of mixture function, i.e.: 
\begin{align}
\label{eq:gradient}
b_{\eta}(\bo{x}) &\eqdef \nabla P_{\eta}(\bo{x})\eqsp,\\
\label{eq:GaP}
P_{\eta}(\bo{x}) &= \sum_{i=1}^K \pi_k\,\varphi^\eta_k(\bo{x}) \quad\mbox{with}\quad \varphi^\eta_k(\bo{x}) \eqdef \exp \left\{-\frac{1}{2}(\bo{x}-\boldsymbol{\mu_k})^T C_k (\bo{x}-\boldsymbol{\mu_k})\right\}\eqsp,
\end{align}
where: 
\begin{enumerate}[-]
\item $K$ is the number of components of the mixture;
\item $\pi_k \in \R^+$ is the relative weight of the $k$-th component, with $\sum_{k=1}^K\pi_k = 1$;
\item  $\boldsymbol{\mu_k} \in \R^2$ is the center of the $k$-th component;
\item $C_k \in \mathcal{S}_2^+$ is the information matrix of the $k$-th component, where $\mathcal{S}_2^+$ is the set of $2\times 2$ symmetric positive definite matrices.
\end{enumerate}
Equations \eqref{eq:target:sde1} and \eqref{eq:gradient} define a potential based continuous time movement model, and follows a popular idea in movement ecology \citep{blackwell:1997,brillinger:et:al:2001a, brillinger:et:al:2001b, brillinger:et:al:2002, preisler:et:al:2004, brillinger:2010,  preisler:et:al:2013, harris:blackwell:2013}. 
In this framework, the movement is supposed to reflect the attractiveness of the environment, which is modeled using a real valued potential function defined on $\mathbb{R}^2$. 
Equation \eqref{eq:GaP} describes the new potential function proposed in this paper. This function is a smooth multimodal surface, having $K$ different components. In this paper, $K$ is supposed to be known one we need to infer $\eta$ from movement data where:
\[
\eta = \{ (\pi_k)_{k=1,\dots,K-1},(\boldsymbol{\mu_k})_{k=1,\dots,K} ,(C_k)_{k=1,\dots,K}\}\in \mathbb{R}^{6K-1}\eqsp.
\]
Potential based continuous time movement model have been widely studied in ecology for the last twenty years. 
The two most popular SDE based models in movement ecology are:
\begin{enumerate}[-]
\item the Brownian motion, with or without drift, corresponding to a constant $b_{\eta}$ (corresponding to a linear potential function), as in \cite{skellam_random_1951} ;
\item the Ornstein-Uhlenbeck process, corresponding to a linear drift function $b_{\eta}$ (corresponding to a linear potential function), as in \cite{blackwell:1997}.
\end{enumerate}
These two models are appealing as their transition density function is available explicitly, and the inference in these case is straightforward. 
However, it contains only a very restricted class of potential functions $P_{\eta}$ (linear, or quadratic) that might be unrealistic in ecological applications (unimodal attractiveness surfaces, for instance). 
 The model introduced in this paper is an extension  of one of those suggested in \cite{preisler:et:al:2013}, it has multiple modes and is flexible enough to accomodate  a wide range of situations. 
 However, the process solution to the SDE \eqref{eq:target:sde1} with the potential function \eqref{eq:GaP} has no explicit transition density so that approximate maximum likelihood procedures such as the methods presented in Section~\ref{sec:inference} may be considered.
\section{Maximum likelihood estimation procedures}
\label{sec:inference}  
The aim of this section is to present four maximum likelihood inference procedures to estimate the parameter $\bs{\theta}=(\eta,\gamma)$ using a set of discrete time observations. 
Among them, the Euler-Maruyama method is the only one used so far in movement ecology, regardless of its potential weaknesses. 
The performance of these methods are evaluated in Section~\ref{sec:application}.

The observations  are given by $G$ independent realizations of the continuous process given by equations \eqref{eq:target:sde1}-\eqref{eq:GaP}  . The set of observations is then the collection of trajectories $(\x^g)_{g=1,\dots,G}$ where each trajectory $\x^g$ is made of $n_g+1$ exact observations  at times $t_0^g=0 < t_1^g < \dots < t_{n^g}^g = T^g$ and  starting at $\bo{x}_0^g$. The probability density of the distribution of $\bo{X}_{\Delta}$ given $\bo{X}_0$ (i.e., the transition density) is denoted by $q^{\Delta}_{\bs{\theta}}$ when the model is parameterized by $\bs{\theta}$: for all bounded measurable function $h$ on $\mathbb{R}^2$:
\[
\mathbb{E}_{\bs{\theta}}\left[h(\bo{X}_{\Delta})\middle| \bo{X}_0\right] = \int h(\bo{y}) q^{\Delta}_{\bs{\theta}}(\bo{X}_0,\bo{y})\rmd \bo{y}\eqsp,
\]
where $\mathbb{E}_{\bs{\theta}}$ is the expectation when the model is driven by $\bs{\theta}$.

As the $G$ trajectories $(\x^g)_{g=1,\dots,G}$ are independent, and by the Markov property of the solution to \eqref{eq:target:sde1}, the likelihood function may be written:
\begin{equation}
\label{eq:target:likelihood}
L(\bs{\theta};{\bf\x})=\prod_{g=1}^G \prod_{i=0}^{n^g-1} q^{\Delta_i^g}_{\bs{\theta}}(\bo{x}_{i}^g,\bo{x}_{i+1}^g)\eqsp,
\end{equation}
where $\Delta_i^g \eqdef t_{i+1}^g-t_{i}^g$. The maximum likelihood estimator of $\bs{\theta}$ is defined as:
\[
\hat{\bs{\theta}} : = \text{argmax}_{\bs{\theta}}\,L(\bs{\theta};{\bf\x})\eqsp.
\]
In the cases of the Brownian motion and Ornstein-Uhlenbeck processes, the function $q^{\Delta}_{\bs{\theta}}$ is the probability density function of a Gaussian random variable with known mean and covariance so that the direct computation of $\hat{\bs{\theta}}$ is an easy task. 
However, in more general settings such as in the model considered in Section~\ref{sec:model}, $q^{\Delta}_{\bs{\theta}}$ is unknown. 
 While many statistical methods have been proposed to compute approximations of $\hat{\bs{\theta}}$, the only statistical method used in movement ecology is the Euler-Maruyama discretization procedure. 
It is known that the quality of this approximation highly depends on the sampling rate which means that this method might therefore not be well suited to real life animal GPS tags.
Therefore, several alternatives are investigated in the remainder of this paper to compute approximations of $\hat{\bs{\theta}}$.
Since all $g$ trajectories are exchangeable, for the sake of clarity, the aforementioned approximation methods are detailed for a given trajectory and references to $g$ are omitted.
\paragraph{Euler-Maruyama method}
The Euler-Maruyama discretization is an order 0 Taylor expansion of the drift function (as the diffusion coefficient is constant no expansion is required). 
The drift is therefore assumed to be constant between two observations. 
For any $0\le i \le n-1$ the target process is approximated by the process $({\bo{X}}^\mathsf{E}_t)_{t_{i}\leq t < t_{i+1}}$, solution to the SDE: ${\bo{X}}^\mathsf{E}_{t_i} = \bo{x_i}$ and, for $t_i\le t\le t_{i+1}$,
\begin{equation*}
\rmd {\bo{X}}^\mathsf{E}_t= b_{\eta}(\bo{x_i})\rmd t + \gamma \rmd \bo{W}_t\eqsp.
\end{equation*}
For all $0\le i \le n-1$, the transition density $q^{\Delta_i}_{\bs{\theta}}(\bo{x_i},\bo{x_{i+1}})$ is approximated by the transition density of the process $({\bo{X}}^\mathsf{E}_t)_{t_i \leq t \leq t_{i+1}}$, denoted by $q^{\Delta_i,\mathsf{E}}_{\bs{\theta}}$,  evaluated for the observations $(\bo{x_{i}},\bo{x_{i+1}})$. 
This is the probability density function of a Gaussian random variable with mean $\boldsymbol{\mu_i}^{\mathsf{E}}$ and variance $\Sigma_i^{\mathsf{E}}$:
\begin{align*}
\boldsymbol{\mu_i}^{\mathsf{E}} &= \bo{x_i} + \Delta_i b_{\eta}(\bo{x_i})\eqsp,\\
\Sigma_i^{\mathsf{E}} &= \gamma^2\Delta_i\eqsp,
\end{align*}
where $I_k$ denotes the $k\times k$ identity matrix.
Hence, the Euler-Maruyama estimate is given by maximizing \eqref{eq:target:likelihood} when $q^{\Delta_i^g}_{\bs{\theta}}$ is replaced by $q^{\Delta_i^g,\mathsf{E}}_{\bs{\theta}}$.
\paragraph{Ozaki method}
The Ozaki method, proposed in \cite{ozaki:1992} and \cite{shoji:1998}, provides a local linearization of the drift term in order to improve the Euler scheme. 
For  any $0\le i \le n-1$ the target process is approximated by the process $({\bo{X}}^\mathsf{O}_t)_{t_{i}\leq t < t_{i+1}}$, solution to the SDE: ${\bo{X}}^\mathsf{O}_{t_i} = \bo{x_i}$ and, for $t_i\le t\le t_{i+1}$,
\begin{equation*}
\rmd {\bo{X}}^\mathsf{O}_t= \left[J_{i,\eta}\left({\bo{X}}^\mathsf{O}_t-\bo{x_i}\right)+ b_{\eta}(\bo{x_i})\right]\rmd t + \gamma \rmd \bo{W}_t\eqsp,
\end{equation*}
where $J_{i,\eta}$ is the $2\times 2$ Jacobian matrix of the drift function $b_{\eta}$ evaluated at $\bo{x_i}$. 
Therefore, the target process is now approximated between each pair of observations by a two-dimensionnal Ornstein-Uhlenbeck process. 
For all $0\le i \le n-1$, the transition density $q^{\Delta_i}_{\bs{\theta}}(\bo{x_i},\bo{x_{i+1}})$ is approximated by the transition density of the process $({\bo{X}}^\mathsf{O}_t)_{t_i \leq t \leq t_{i+1}}$, denoted by $q^{\Delta_i,\mathsf{O}}_{\bs{\theta}}$,  evaluated for the observations $(\bo{x_{i}},\bo{x_{i+1}})$. 
If the potential function  is such that, for all $0\le i\le n-1$ and all $\bs{\theta}$, $J_{i,\bs{\theta}}$ is nonsingular and symmetric (as it is the case of the potential introduced in Section~\ref{sec:model}), following \cite{shojiOzaki1998},  this transition density is the probability density function of a Gaussian  random variable with mean $\boldsymbol{\mu_i}^{\mathsf{O}}$ and variance $\Sigma_i^{\mathsf{O}}$:
\begin{align*}
\boldsymbol{\mu_i}^{\mathsf{O}} &= \bo{x_i} + (\exp\left(J_{i,\eta} \Delta_i \right) - I_2) (J_{i,\eta})\inv b_{\eta}(\bo{x_i})\eqsp,\\
\textbf{vec}(\Sigma_i^{\mathsf{O}}) &= (J_{i, \eta} \oplus J_{i, \eta} )\inv \left(e^{(J_{i,\eta} \oplus J_{i,\eta})\Delta_i} - I_4\right)\text{vec}\left(\gamma^2 I_2\right)\eqsp,
\end{align*}
where \textbf{vec} is the stack operator and $\oplus$ is the Kronecker sum.
The local linearisation estimate is then given by maximizing \eqref{eq:target:likelihood} when $q^{\Delta_i^g}_{\bs{\theta}}$ is replaced by $q^{\Delta_i^g,\mathsf{O}}_{\bs{\theta}}$.

\paragraph{Kessler method}
This procedure, proposed in \cite{kessler:1997}, generalizes the Euler-Maruyama method using a higher order approximation of the transition density. 
Following the same steps as Euler or Ozaki discretizations,  \cite{kessler:1997} proposed to consider a Gaussian approximation of the transition density between consecutive observations. 
The procedure aims at replacing $q^{\Delta_i}_{\bs{\theta}}(\x_{i},\x_{i+1})$ by a Gaussian probability density function with mean and variance given by the exact mean $\boldsymbol{\mu_i}^{\mathsf{target}}$ and variance $\Sigma_i^{\mathsf{target}}$ of the target process at time $t_{i+1}$ conditionally on the value of the target process at time $t_i$. 
This leads to an approximation with the same first two moments as the target transition density. 
This approximation can be written similarly to the previous discretization schemes by defining, for  any $0\le i \le n-1$, the process $({\bo{X}}^\mathsf{target}_t)_{t_{i}\leq t < t_{i+1}}$ as the solution to the SDE: ${\bo{X}}^\mathsf{target}_{t_i} = \bo{x}_{i}$ and, for $t_i\le t\le t_{i+1}$,
\begin{equation*}
\rmd {\bo{X}}^\mathsf{target}_t= \frac{\boldsymbol{\mu_i}^{\mathsf{target}} - \bo{x_i} }{\Delta_i}\rmd t + \left(\frac{\Sigma_i^{\mathsf{target}}}{\Delta_i}\right)^{1/2} \rmd \bo{W}_t\eqsp.
\end{equation*}
However, the transition density of the target process being unknown, $\boldsymbol{\mu_i}^{\mathsf{target}}$ and $\Sigma_i^{\mathsf{target}}$ are not available in closed form.
 Nonetheless, Taylor expansions of these two conditionnal moments are available (see  \cite{florens1989approximate}, or \cite{uchida:yoshida:2012} for the multidimensional case).
  \cite{kessler:1997} proposed to replace $\mu_i^{\mathsf{target}}$ and $\Sigma_i^{\mathsf{target}}$ by these approximations. 
For all $0\le i \le n-1$, the transition density $q^{\Delta_i}_{\bs{\theta}}(\bo{x_i},\bo{x_{i+1}})$ is approximated by the transition density of the process $({\bo{X}}^\mathsf{K}_t)_{t_i \leq t \leq t_{i+1}}$ (approximated with these expansions). 
These expansions can be computed directly from the drift and diffusion functions and their partial derivatives. 
The order of the expansion of the true conditional moments is left to the user and the performance of the estimator highly depends on this parameter, order one being the Euler-Maruyama method. 
In this paper, it was performed up to the second order. 
In this case, and for the process defined by \eqref{eq:target:sde1} and \eqref{eq:gradient}, the function $q^{\Delta_i,\mathsf{K}}_{\bs{\theta}}(\bo{x_i},\cdot)$ is the probability density function of a Gaussian random variable with mean $\boldsymbol{\mu_i}^{\mathsf{K}}$ and variance $\Sigma_i^{\mathsf{K}}$:
\begin{align}
\boldsymbol{\mu_i}^{\mathsf{K}} &= \bo{x_i} + \Delta_i b_{\eta}(\bo{x_i})\eqsp,\nonumber\\
\Sigma_i^{\mathsf{K}} &= \gamma^2\Delta_i(I_2+\Delta_i J_{i,\eta})\eqsp.\label{eq:kes:var}
\end{align}
Note that $\boldsymbol{\mu_i}^{\mathsf{K}}=\boldsymbol{\mu_i}^{\mathsf{E}}$. 
For numerical stability, $(\Sigma_i^{\mathsf{K}})^{-1}$ and its determinant can be replaced by a Taylor expansion, the associated contrast function to optimize in multidimensional cases is given in \cite{uchida:yoshida:2012}.
 A drawback of this method is the fact that \eqref{eq:kes:var} does not necessarily define a positive semi-definite matrix. For instance, if $\text{Tr}(J_{i,\eta}) < -2/\Delta_i$, then $\text{Tr}(\Sigma_i^{\mathsf{K}})<0$, which is likely to occur when $\Delta_i$ is large. In the following applications, whenever $\Sigma_i^{\mathsf{K}}$ is not positive definite, the associated observation is thrown out from the computation of the likelihood, as proposed in \cite{iacus2009simulation}. Kessler estimate is then given by maximizing \eqref{eq:target:likelihood} when $q^{\Delta_i^g}_{\bs{\theta}}$ is replaced by $q^{\Delta_i^g,\mathsf{K}}_{\bs{\theta}}$.
\paragraph{Exact Algorithm based Monte Carlo EM (EA MCEM) method}
The EA MCEM proposed by \cite{beskos:papaspiliopoulos:roberts:fearnhead:2006} does not use  Gaussian approximations of the transition densities to compute the maximum likelihood estimator. 
Details of the method are given in \cite{beskos:papaspiliopoulos:roberts:fearnhead:2006} and the important results for the experiments are given online as a supplementary material.
 
Applying the exact algorithm to estimate $\bs{\theta}$ requires the target SDE to be reducible to a unit diffusion using the Lamperti transform which is obtained in our case by setting $(\bo{Y}_t:= \gamma\inv \bo{X}_t)_{0\le t\le T}$. Then,
\begin{equation}
\label{eq:target:sdeY}
\rmd \bo{Y}_t=\alpha_{\bs{\theta}}(\bo{Y}_t)\rmd t + \rmd \bo{W}_t\quad\mbox{where}\quad \alpha_{\bs{\theta}}(\cdot) \eqdef \gamma\inv b_\eta(\gamma\,\cdot) = \gamma\inv\nabla P_{\eta}(\gamma\,\cdot)\eqsp.
\end{equation}
The EA MCEM approach used in this work relies on the following assumptions.
\begin{enumerate}[-]
\item \textit{Conservative assumption}: for all $\bs{\theta}=(\eta,\gamma)\in\mathbb{R}^{d+1}$, there exists $H_{\bs{\theta}}: \mathbb{R}^2 \mapsto \R$ such that for all $\bo{x} \in \mathbb{R}^2$,
\begin{equation}
\alpha_{\bs{\theta}}(\bo{x}) = \nabla H_{\bs{\theta}} (\bo{x})\eqsp.\label{eq:cons:cond}
\end{equation}
\item \textit{Boundedness condition}: for all $\bs{\theta}\in\mathbb{R}^{d+1}$, there exist $m_{\bs{\theta}}, M_{\bs{\theta}}$, such that for all $\bo{x} \in \mathbb{R}^2$,
\begin{equation}
m_{\bs{\theta}} \leq  \|\alpha_{\bs{\theta}}(\bo{x})\|^2 + \Delta H_{\bs{\theta}}(\bo{x} ) \leq M_{\bs{\theta}}\eqsp,\label{eq:bound:cond}
\end{equation}
where $\Delta$ is the Laplace operator:
\begin{equation*}
\Delta H_{\bs{\theta}} : \bo{x}\mapsto \frac{\partial\alpha_{\bs{\theta},1}}{\partial x_1}(\bo{x}) + \frac{\partial\alpha_{\bs{\theta},2}}{\partial x_2}(\bo{x})\eqsp.
\end{equation*}
\end{enumerate}
Both conditions are somehow restrictive in general, however, the first one turns out to be satisfied automatically for potential based models studied in ecology, and described in Section \ref{sec:model}. 
Given the model described by equations \eqref{eq:target:sde1}-\eqref{eq:GaP}, $H_{\bs{\theta}}$ is given as: 
\begin{equation*}
H_{\bs{\theta}}: \bo{x}\mapsto \sum_{k=1}^K \pi_k\varphi^{\eta}_k(\gamma \bo{x})/\gamma\eqsp,
\end{equation*}
with $\bs{\theta}=(\eta, \gamma)$. 
The second condition can be relaxed at the cost of additional computations \citep{beskos:papaspiliopoulos:roberts:fearnhead:2006}, but is satisfied by the potential function of equation~\eqref{eq:GaP}.

This method offers the advantage of avoiding approximation errors due to time discretization of the SDE. 
The only error comes from the Monte Carlo simulations used to approximate the expectation in the E-step (see details in the supplementary material).
In theory, if the conditions \eqref{eq:cons:cond}  and \eqref{eq:bound:cond} are satisfied, it is always possible to perform these simulations. 
However, the simulation algorithm is based on rejection sampling and the acceptance probability is equal to:
\begin{equation}
\exp\left\{-\frac{1}{2}(M_{\bs{\theta}}-m_{\bs{\theta}})\int_{t_i}^{t_{i+1}}\phi_{\bs{\theta}}(\bo{Y}_s)\rmd s\right\}\eqsp,
\label{eq:prob:acc}
\end{equation}
where $m_{\bs{\theta}}$ and $M_{\bs{\theta}}$ are defined as in \eqref{eq:bound:cond}, and  $\phi_{\bs{\theta}}(\bo{x}):= (\parallel \alpha_{\bs{\theta}}(\bo{x}) \parallel^2 + \Delta H_{\bs{\theta}}(\bo{x}) - m_{\bs{\theta}})/2 \geq 0$.
Then, the integral in \eqref{eq:prob:acc} increases with the time step, and the acceptance probability decreases to 0. Therefore, in practice, this method also depend on the sampling step to be used in a reasonnable time.

In order to compare the performance of the three methods, an approximated likelihood criterion is used. The loglikelihood of two given estimates $\hat{\bs{\theta}}_1$ and $\hat{\bs{\theta}}_2$ cannot be computed exactly. 
However, in the case where \eqref{eq:cons:cond} and  \eqref{eq:bound:cond} hold, and following \cite{beskos:papaspiliopoulos:roberts::2009}, an unbiased estimator $\ell_{\mathsf{EA}}$ of the loglikelihood based on Monte Carlo simulations may be computed.
This estimator has no intrinsic bias as the error only comes from the Monte Carlo procedure. In the following, this criterion allowed us to compare the quality of the different estimations resulting from the four different methods.
\section{Application}
\label{sec:application}
\subsection{Simulation study}
\label{sec:sim:study}
\begin{figure}[p]
 \centering
\begin{subfigure}[t]{0.49\textwidth}
\caption{Hidden potential map}
\includegraphics[width=0.99\textwidth]{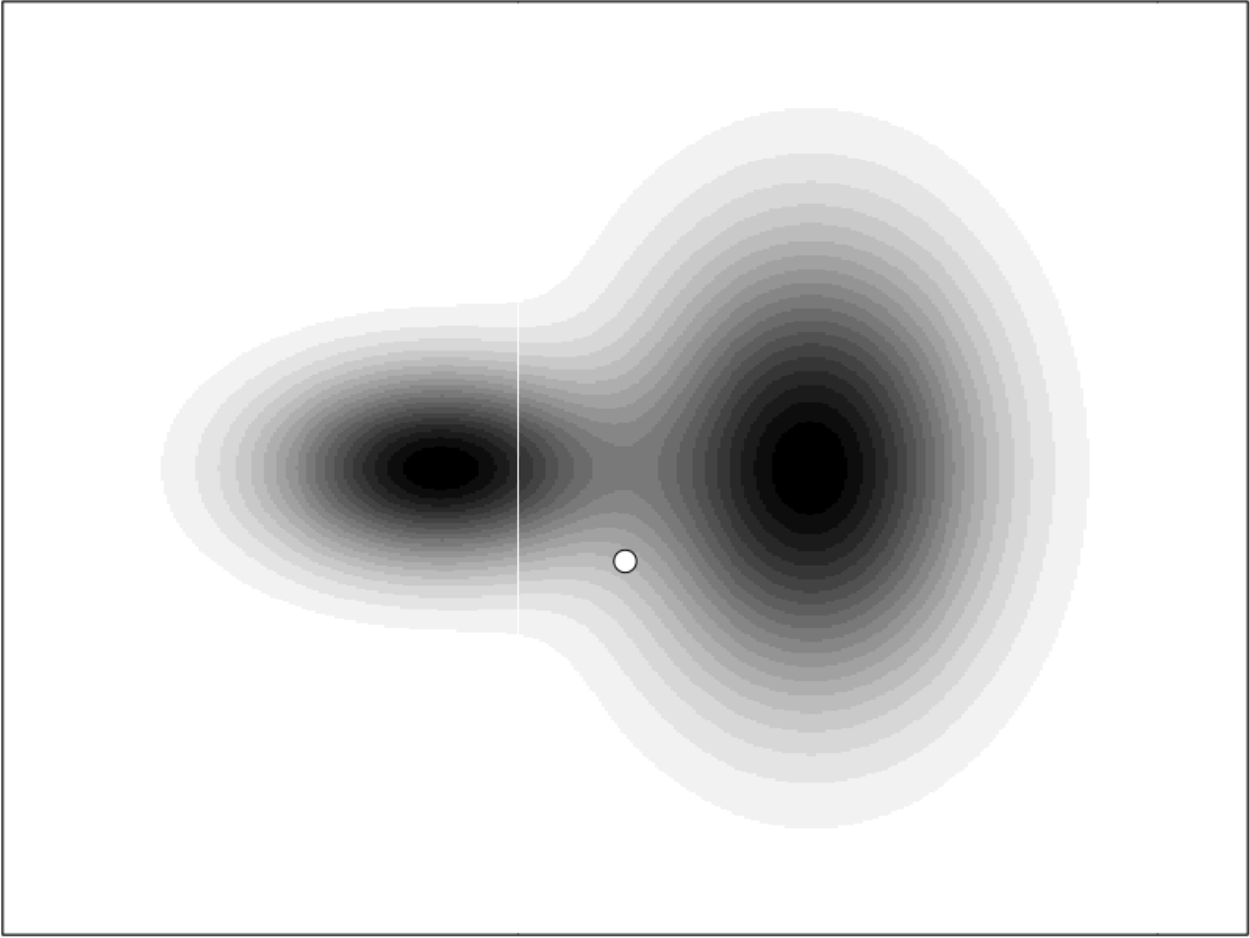}
\label{fig:true:map}
\end{subfigure}
\begin{subfigure}[t]{0.49\textwidth}
\caption{$\Delta =0.1$}
\includegraphics[width=0.99\textwidth]{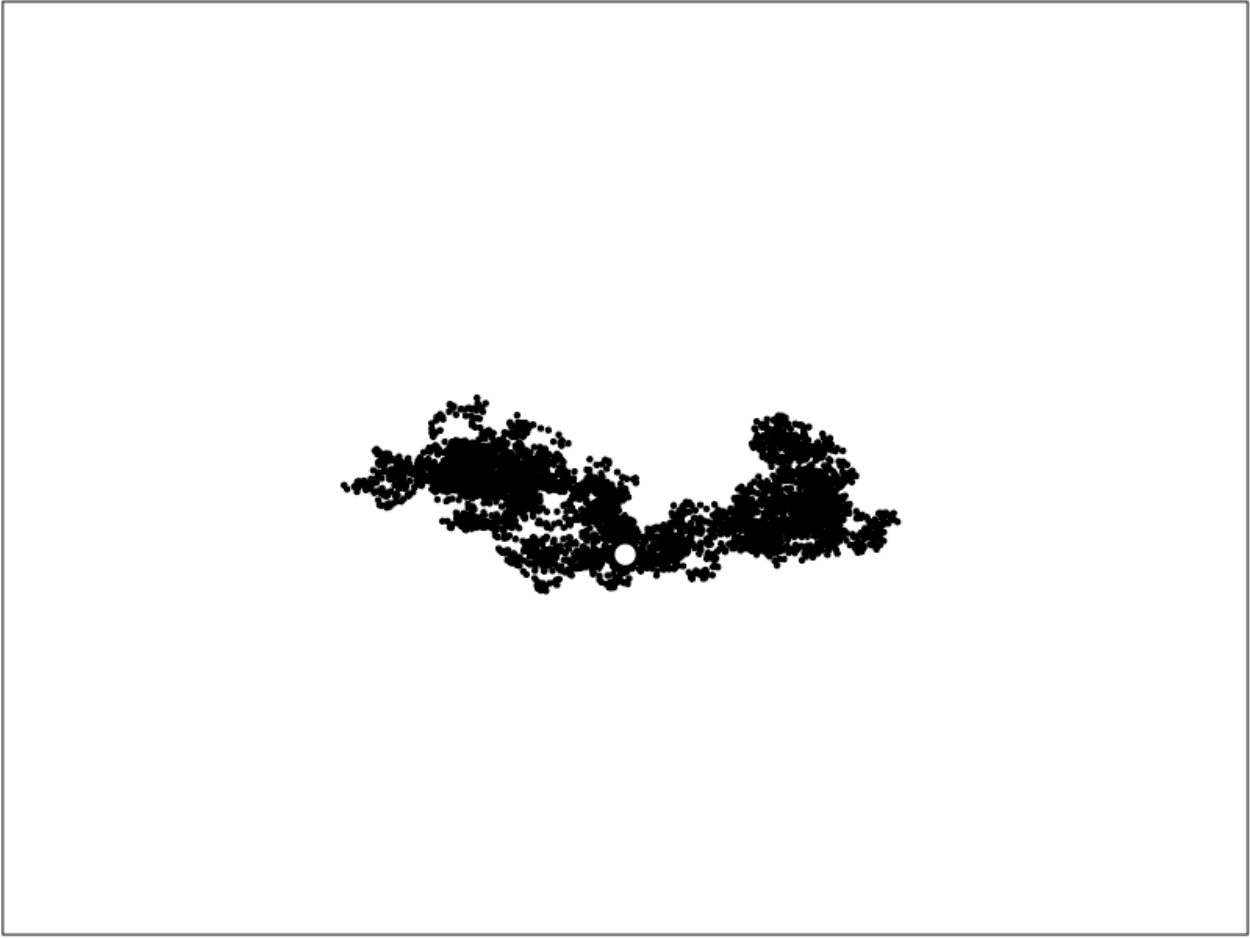}
\end{subfigure}\\
\begin{subfigure}{0.49\textwidth}
\includegraphics[width=0.99\textwidth]{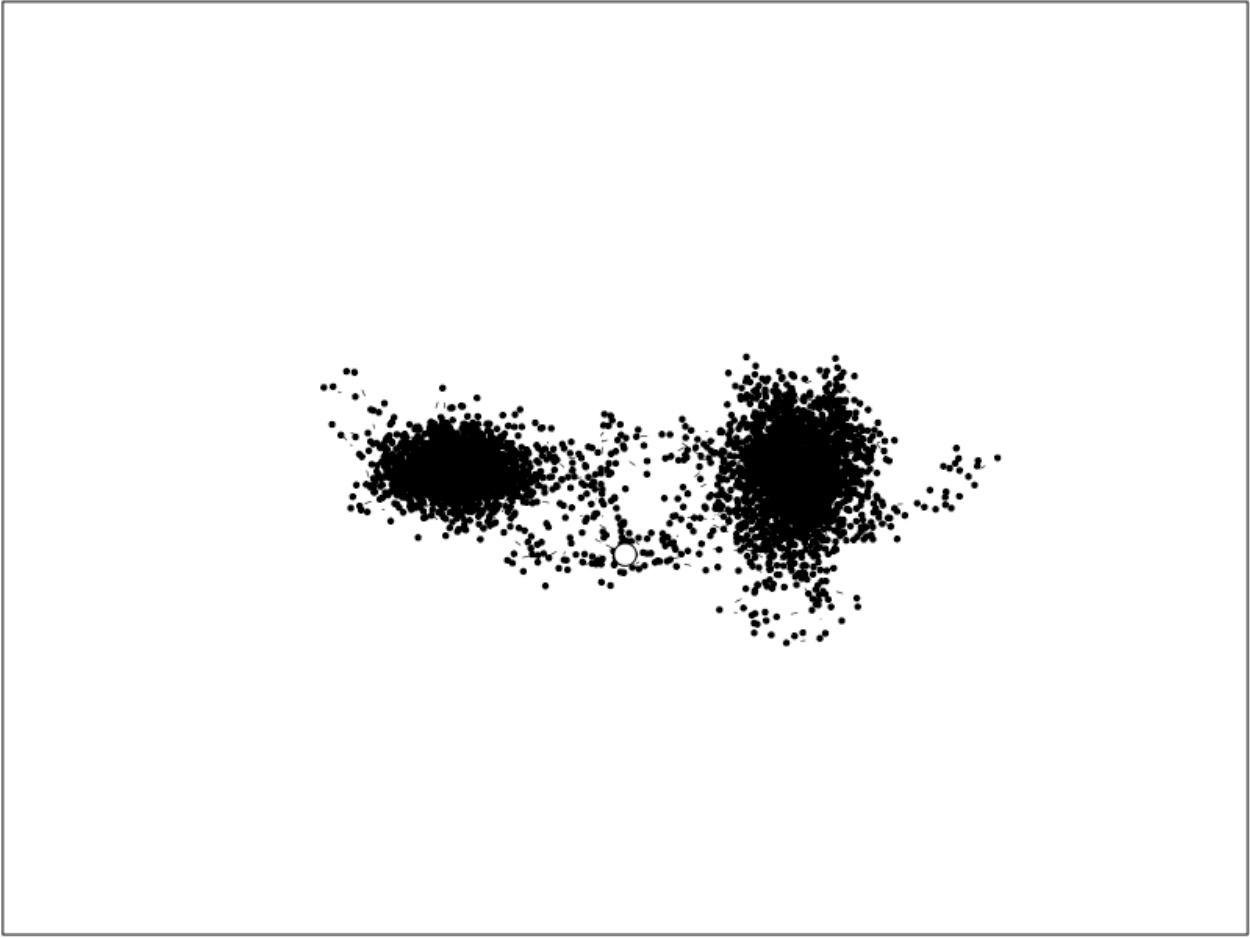}
\caption{$\Delta =1$}
\end{subfigure}
\begin{subfigure}{0.49\textwidth}
\includegraphics[width=0.99\textwidth]{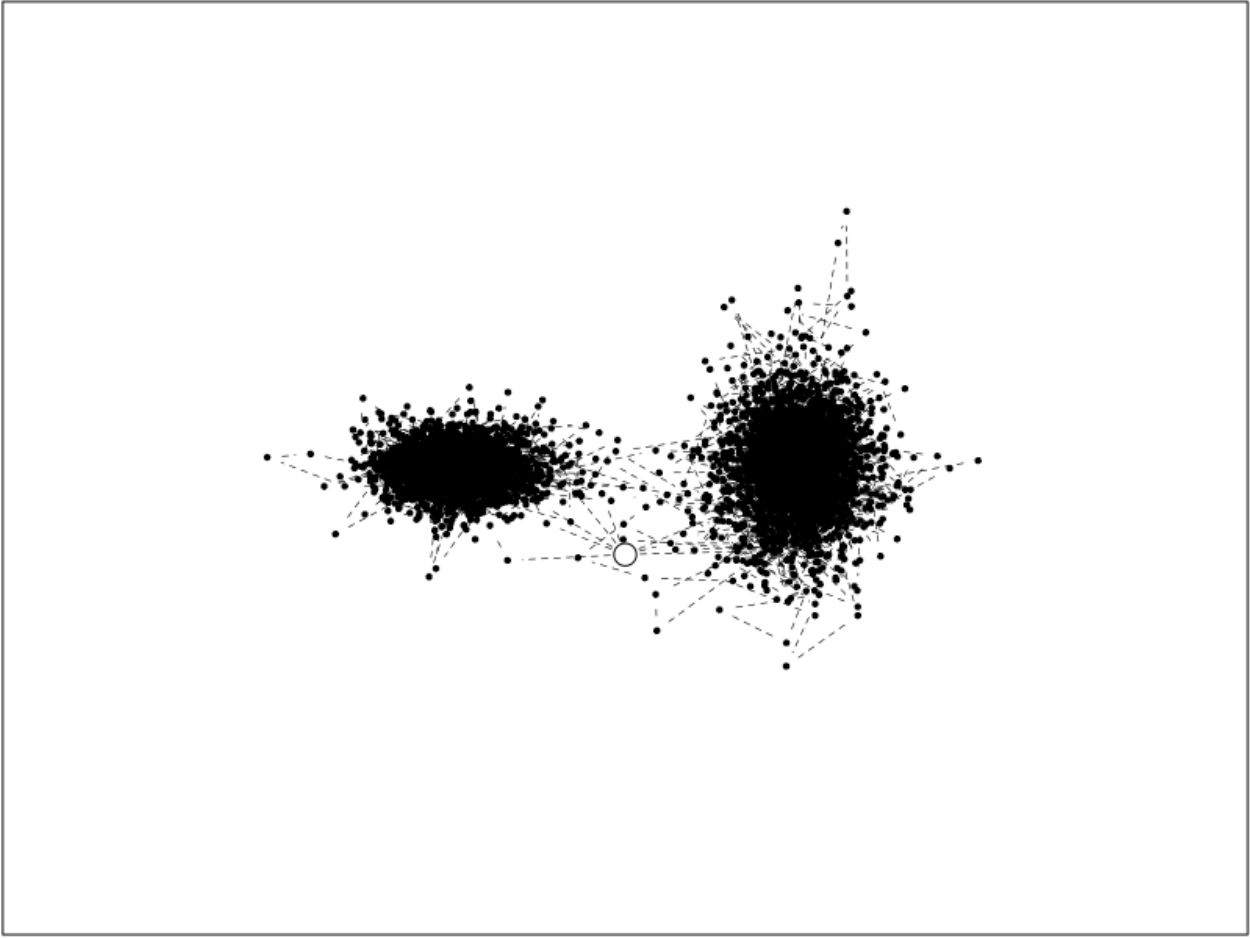}
\caption{$\Delta =10$}
\end{subfigure}
\caption{Simulation of the process solution to \eqref{eq:target:sde1} with a potential as in equation \eqref{eq:GaP}. 
The potential map driving the movement is shown in Figure \ref{fig:true:map},  dark zones present high potential whereas white zones have low potential. On the three other panels, ten simulated trajectories with three different sampling schemes are given. The white dot represents the starting point $\mathsf{x}_0$ of each trajectory.}
\label{fig:sim:data}
\end{figure}

The performance of the four estimation methods presented in this paper is evaluated using different sampling schemes. 
The hidden potential map has $K=2$ components and is represented in Figure \ref{fig:true:map}. 
All simulated trajectories start from the same $\bo{x_0}$ and are sampled from the true distribution of the process $(X_t)_{t\geq 0}$ with parameters $\eta^*$ and $\gamma^*$ using the acceptance rejection procedure proposed in \cite{beskos:roberts:2005}. 
$G = 10$ independent trajectories are simulated with $n_g = 500$ recorded points. 
Three different sampling time steps are considered, respectively with high frequency ($\Delta = 0.1$), intermediate frequency ($\Delta = 1$) and low frequency ($\Delta = 10$). 
In an actual sampling design, this would correspond to three different GPS tags settings with a limited number of emissions. For each method and in each scenario, this procedure is repeated independently to produce 30 different datasets.

\begin{figure}[p]
 \centering
 \includegraphics[width=0.49\textwidth]{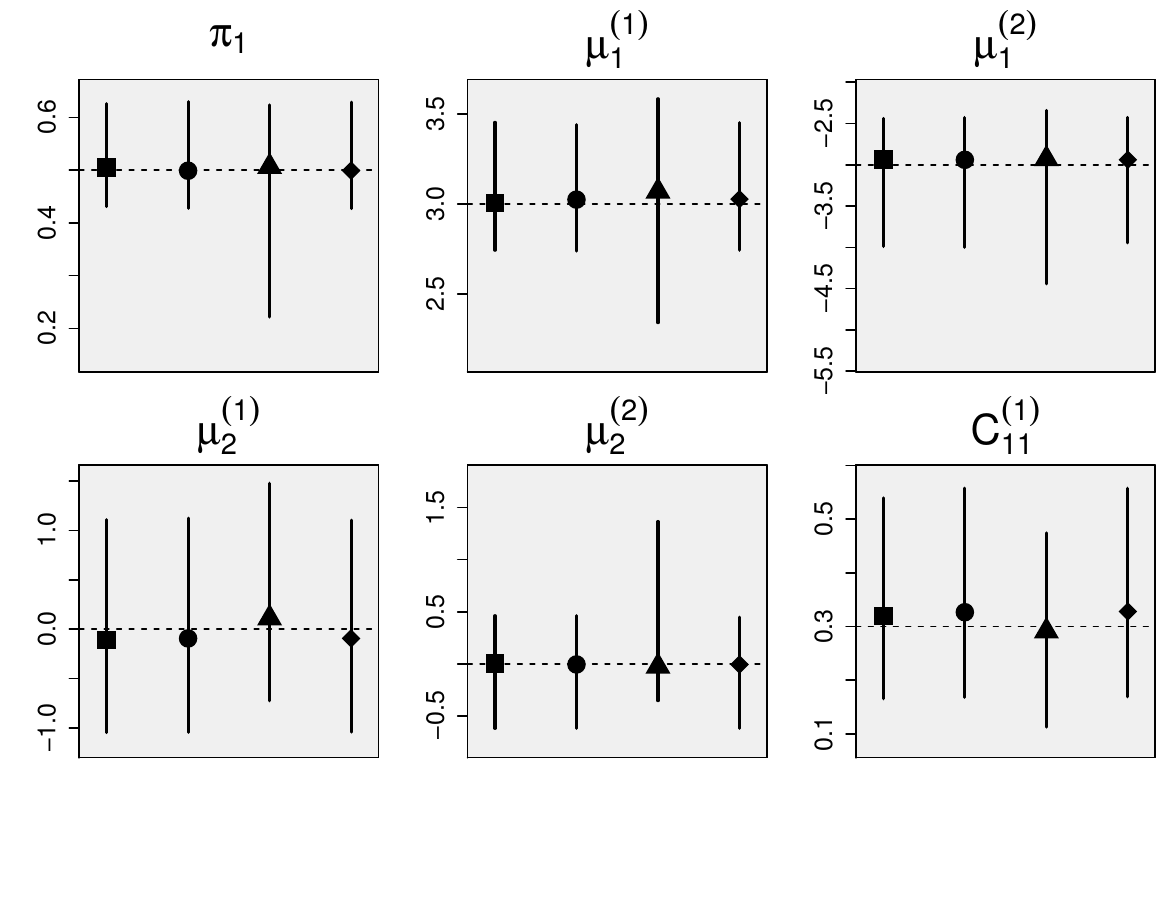}\\
  \vspace{-1.5\baselineskip}
 \includegraphics[width=0.49\textwidth]{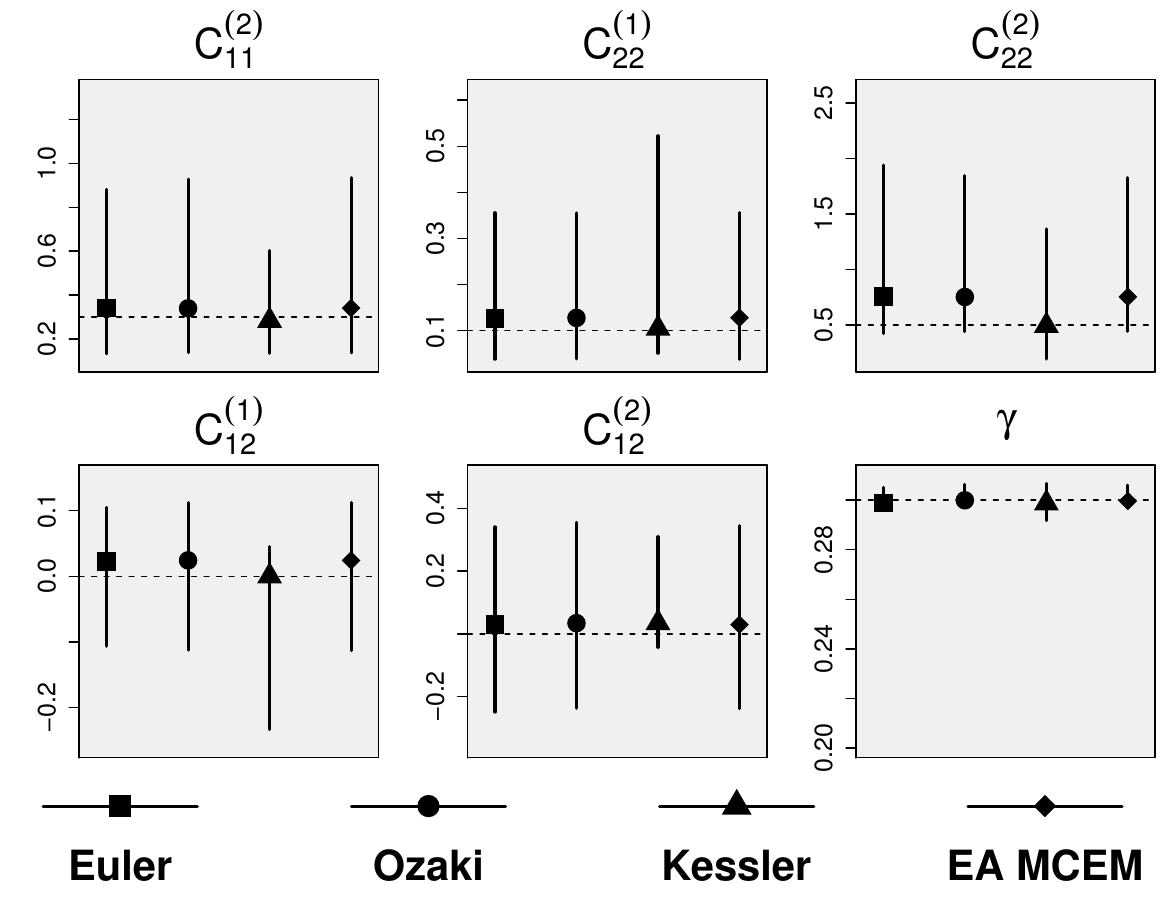}
 \caption{Estimation of the parameters in the case where $\Delta = 0.1$. The dot represents the median, the whiskers provide the 95\% range of estimations. The dashed line represents the true value.}
  \label{fig:est:params:hn01}
 \end{figure}

\begin{figure}[p]
 \centering
 \includegraphics[width=0.49\textwidth]{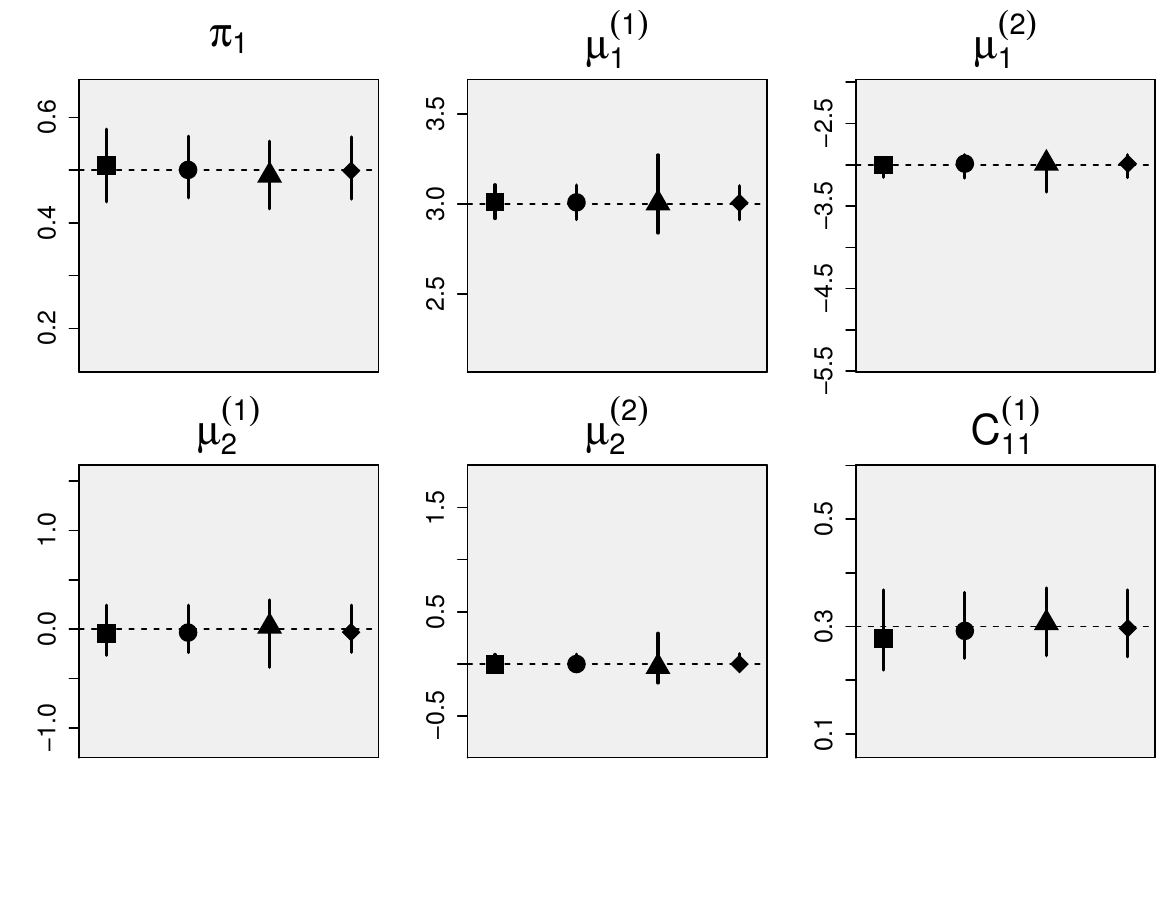}\\
  \vspace{-1.5\baselineskip}
 \includegraphics[width=0.49\textwidth]{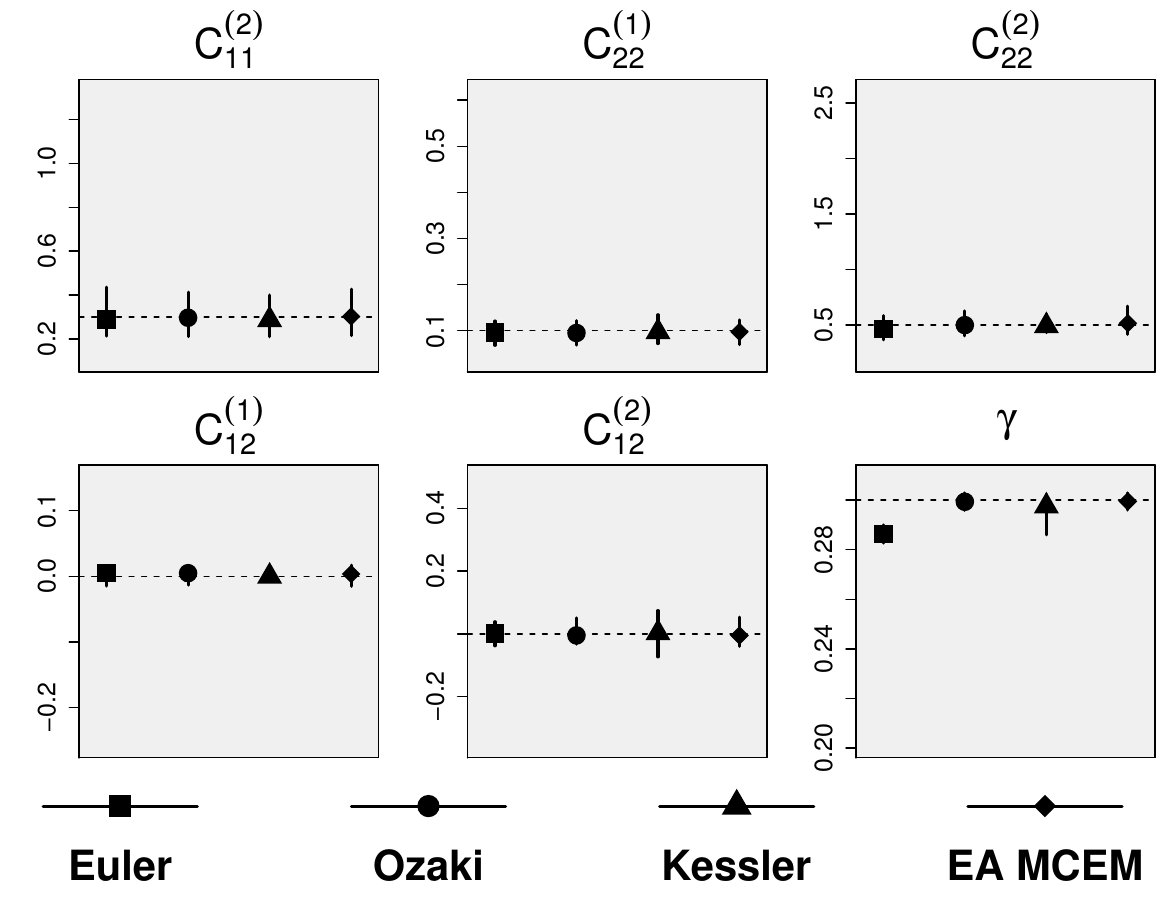}
 \caption{Estimation of the parameters in the case where $\Delta = 1$. The dot represents the median, the whiskers provide the 95\% range of estimations. The dashed line represents the true value.}
  \label{fig:est:params:hn1}
 \end{figure}

\begin{figure}[p]
 \centering
 \includegraphics[width=0.49\textwidth]{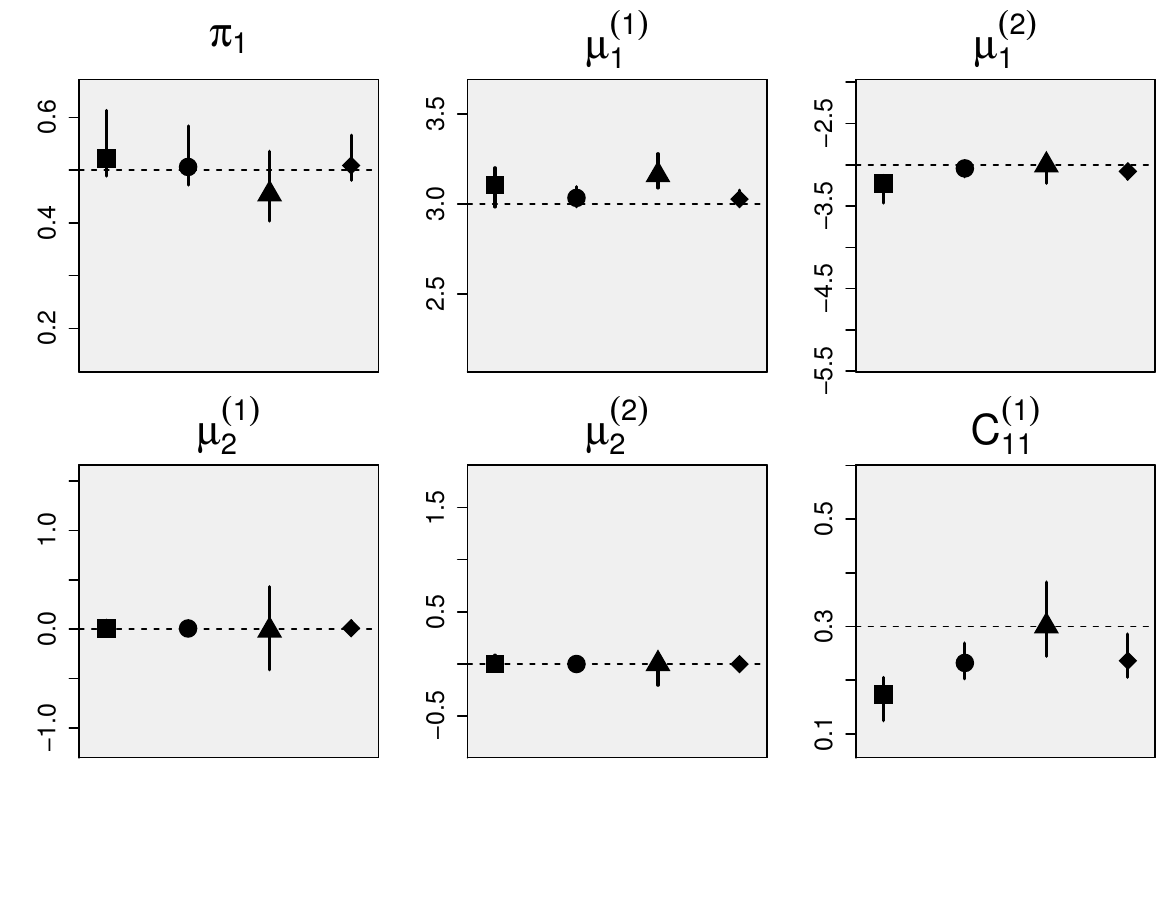}\\
 \vspace{-1.5\baselineskip}
 \includegraphics[width=0.49\textwidth]{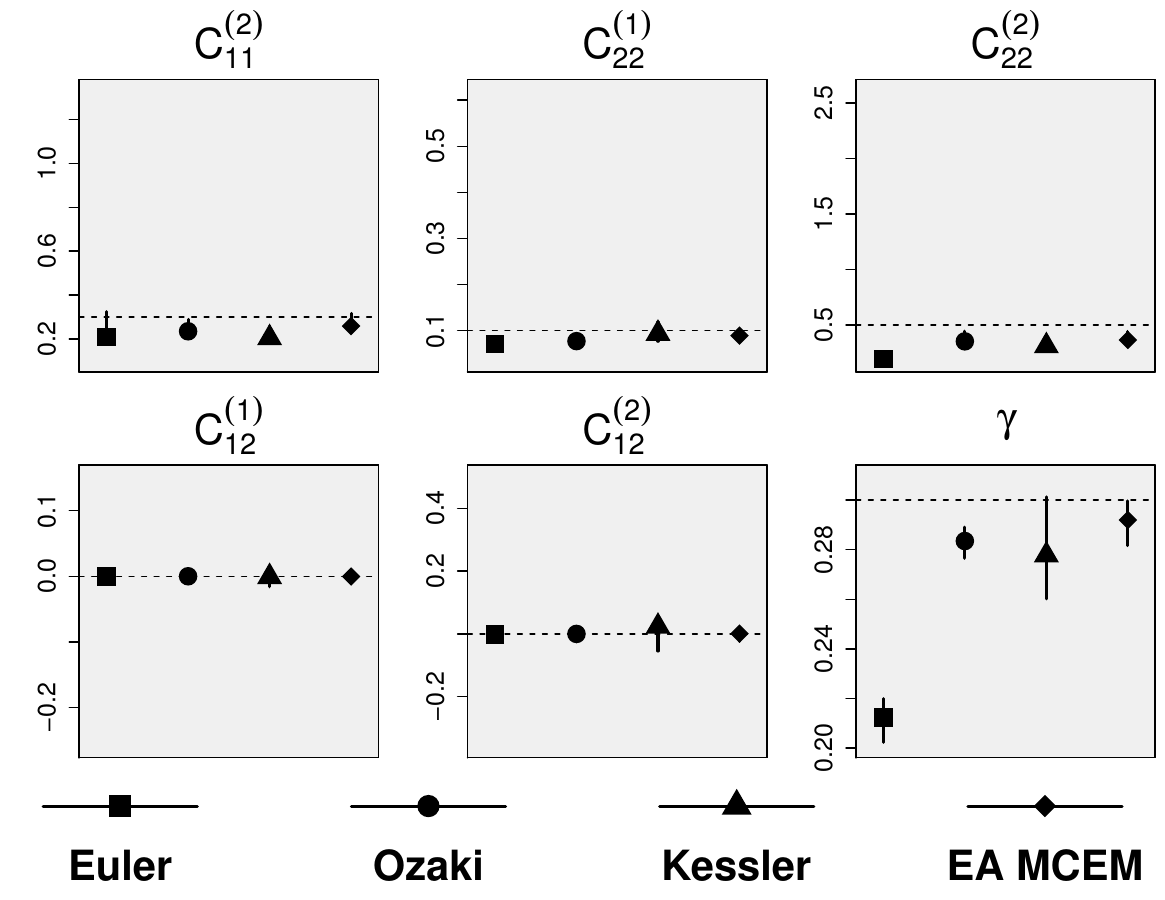}
 \caption{Estimation of the parameters in the case where $\Delta = 10$. The dot represents the median, the whiskers provide the 95\% range of estimations. The dashed line represents the true value.}
  \label{fig:est:params:hn10}
 \end{figure}

The optimization problem associated with the M-step of the EM algorithm  is performed using the CMA-ES algorithm proposed by \cite{hansen_evolutionary_2001}.  For each estimation method, 30 initial parameters $\bs{\theta}_0$ are used and the parameter estimate with the maximum objective function is considered as the estimator $\hat{\bs{\theta}}$.
Figures \ref{fig:est:params:hn01}, \ref{fig:est:params:hn1} and \ref{fig:est:params:hn10} display the parameter estimate for the three sampling cases.
The best estimate for all parameters and techniques is obtained with the intermediate frequency sampling rate. 
As illustrated by Figure \ref{fig:sim:data}, this frequency allows a good exploration of the map by the process with 500 observations and all methods provide reliable estimates.
When $\Delta=0.1$ with only 500 observations per trajectory, the time horizon $50$  might be insufficient for a good parameter estimation using this model.  
On the other hand, when $\Delta = 10$,  there is a strong bias in the estimation of the shape parameters $C^{k}, k=1,2$, and the diffusion coefficient for all methods. 

The Euler method  performs much worse  than the other methods when $\Delta = 1$ and $\Delta = 10$  whereas the three other methods are more robust to the sampling frequency. 
The Ozaki method provides similar results to the EA MCEM method.
The Kessler estimation method provides results which may differ for some parameters (for instance $\boldsymbol{\mu_1}^{(1)}$ or $C_{11}^{(1)}$ in Figure \ref{fig:est:params:hn10}) as about 50\% of the observations lead to a non positive semi-definite matrix in equation \eqref{eq:kes:var} and are therefore dropped. This is a major drawback of this method as it  induces a supplementary bias. 
Finally, note that the EA MCEM method has a much larger computation cost than all other procedures but it is less sensitive to the starting point of the algorithm.
The mean acceptance rate decreases from 99.6\% ($\Delta=0.1$) to 97.6 ($\Delta=1$) and 67.5\% ($\Delta=10$) respectively.

\begin{figure}[p]
 \centering
 \begin{tabular}{ccc}
 \includegraphics[width=0.33\textwidth]{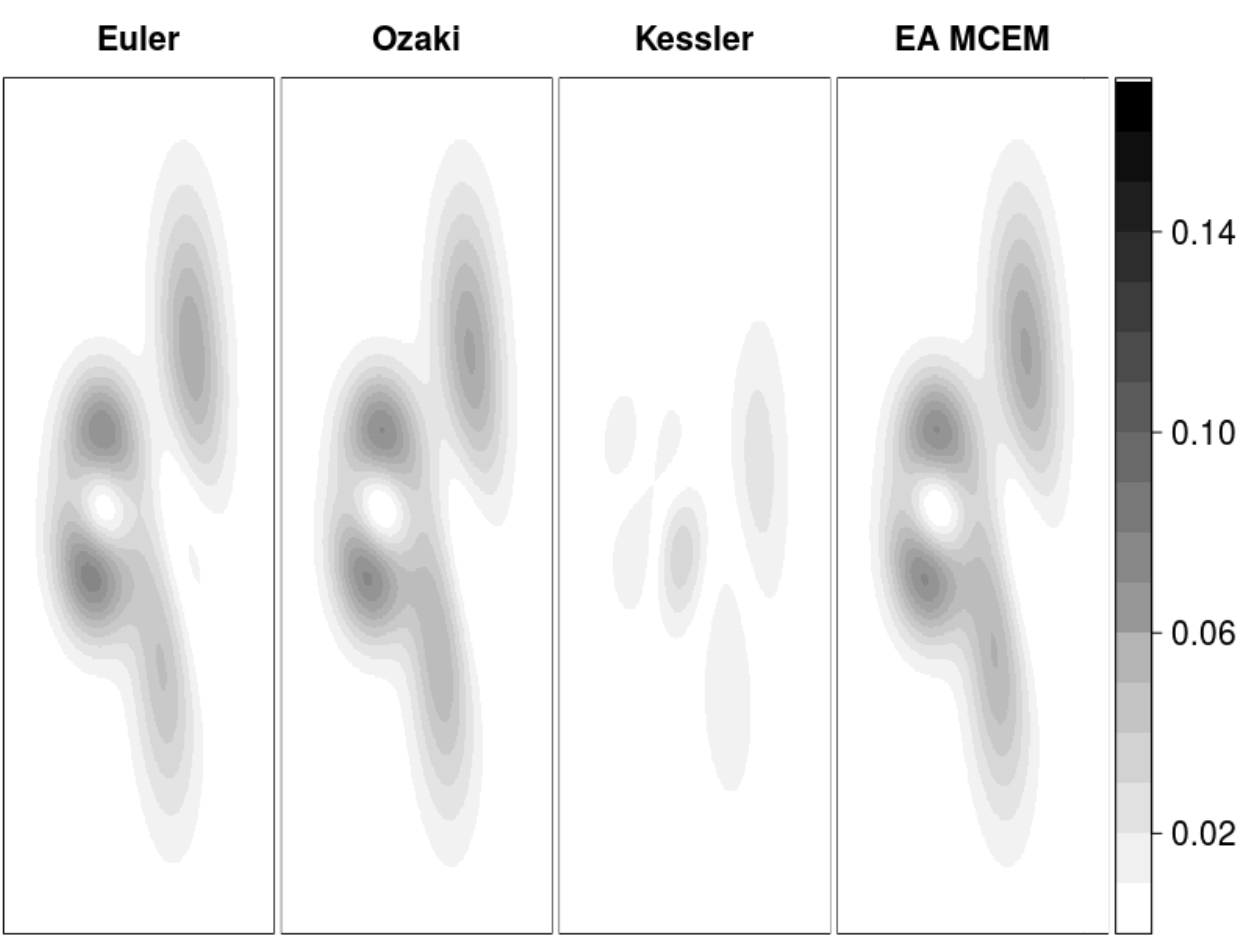}
 \includegraphics[width=0.33\textwidth]{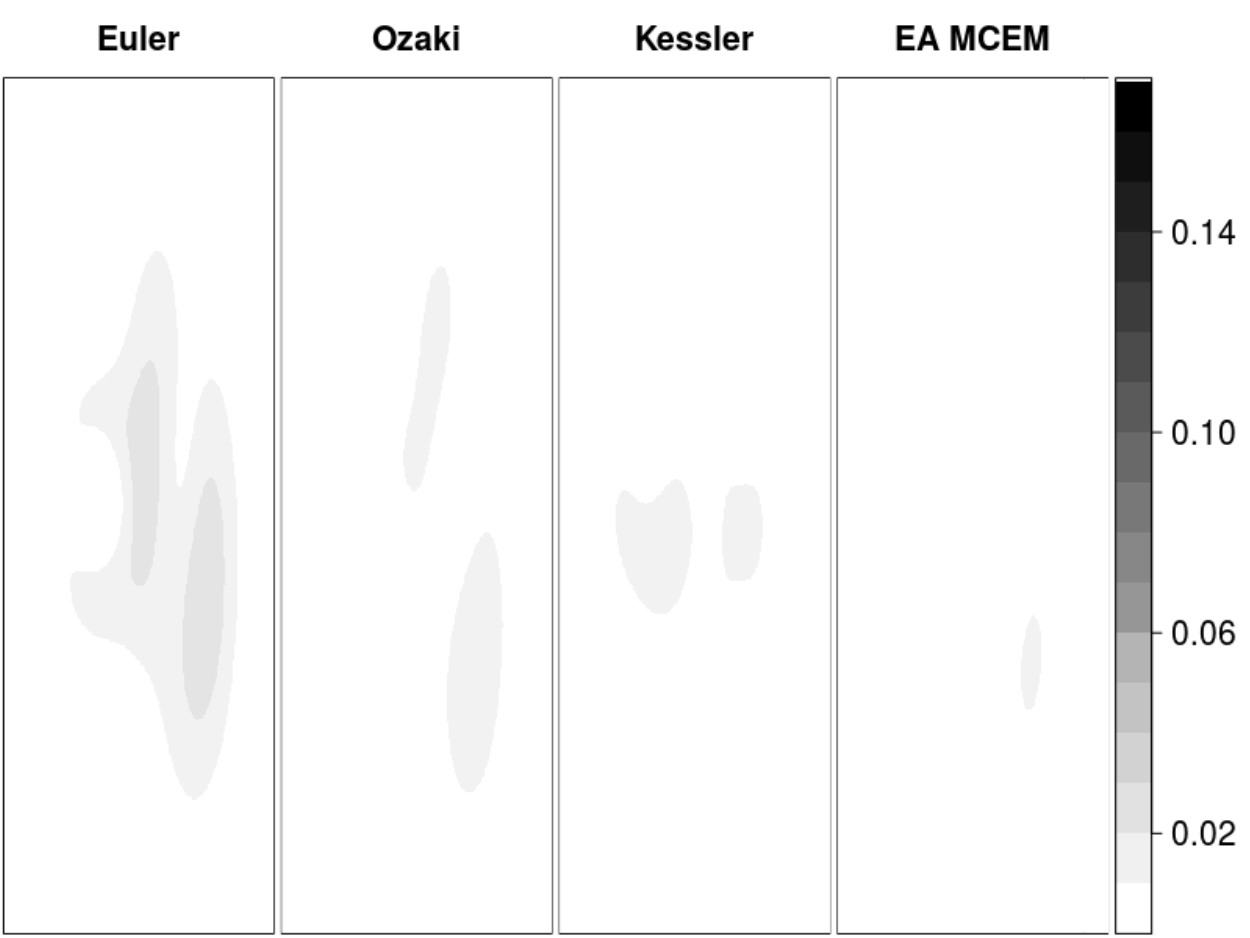}
 \includegraphics[width=0.33\textwidth]{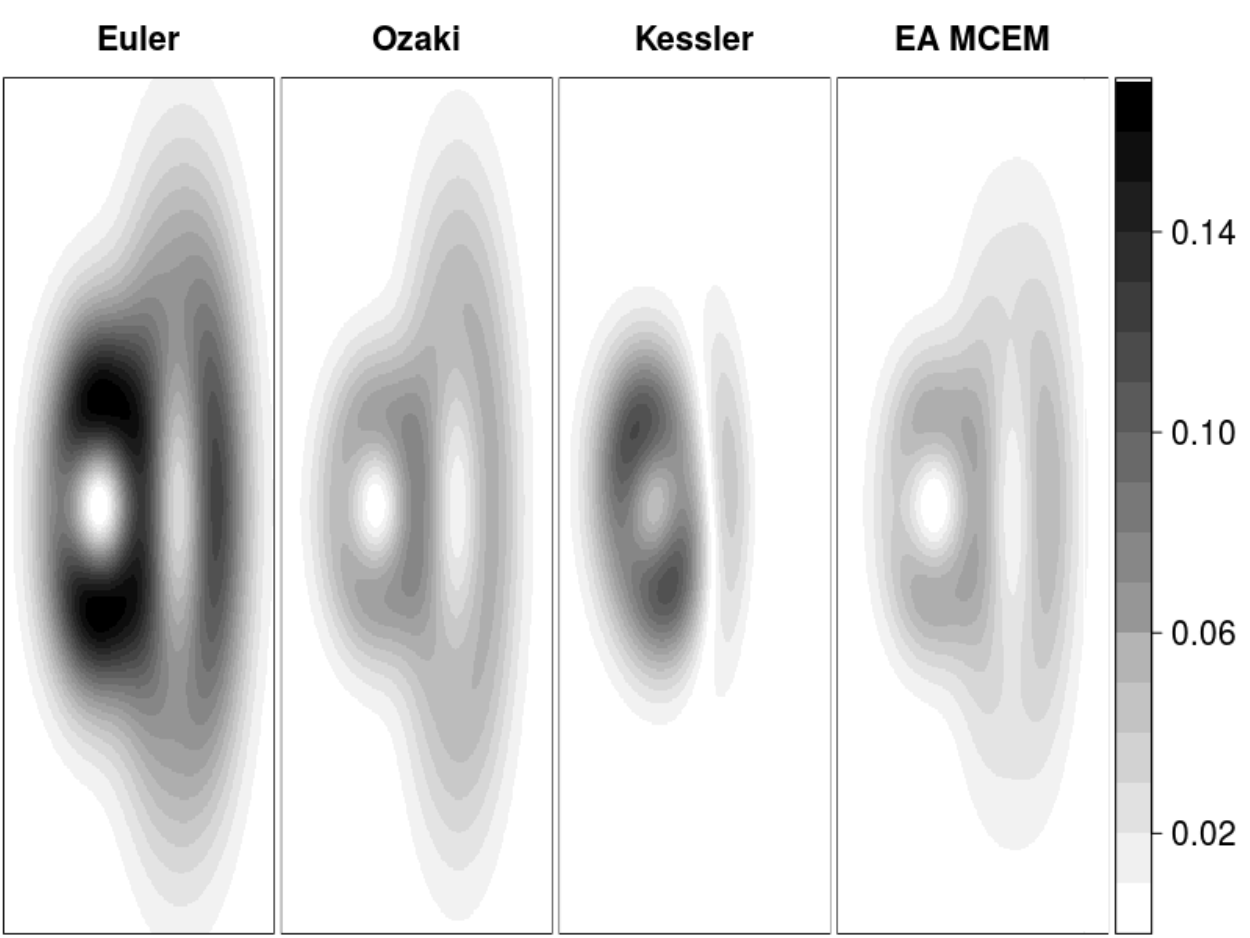}
 \end{tabular}
 \caption{Error between the median map estimate and the true map. Results are presented for each method with $\Delta=0.1$ (left), 1 (center) and 10 (right).}
  \label{fig:map:error}
 \end{figure}

Figure~\ref{fig:map:error} displays the absolute estimation errors $\vert P_{\hat{\eta}}(\bo{x}) - P_{\eta^*}(\bo{x})\vert $ for each method. 
The map is produced using the median of all estimated values as a set of parameters. 
The Euler method  performs much worse  than the other methods when $\Delta = 1$ and $\Delta = 10$. 

\begin{figure}[p]
 \centering
 \begin{tabular}{ccc}
 \includegraphics[width=0.33\textwidth]{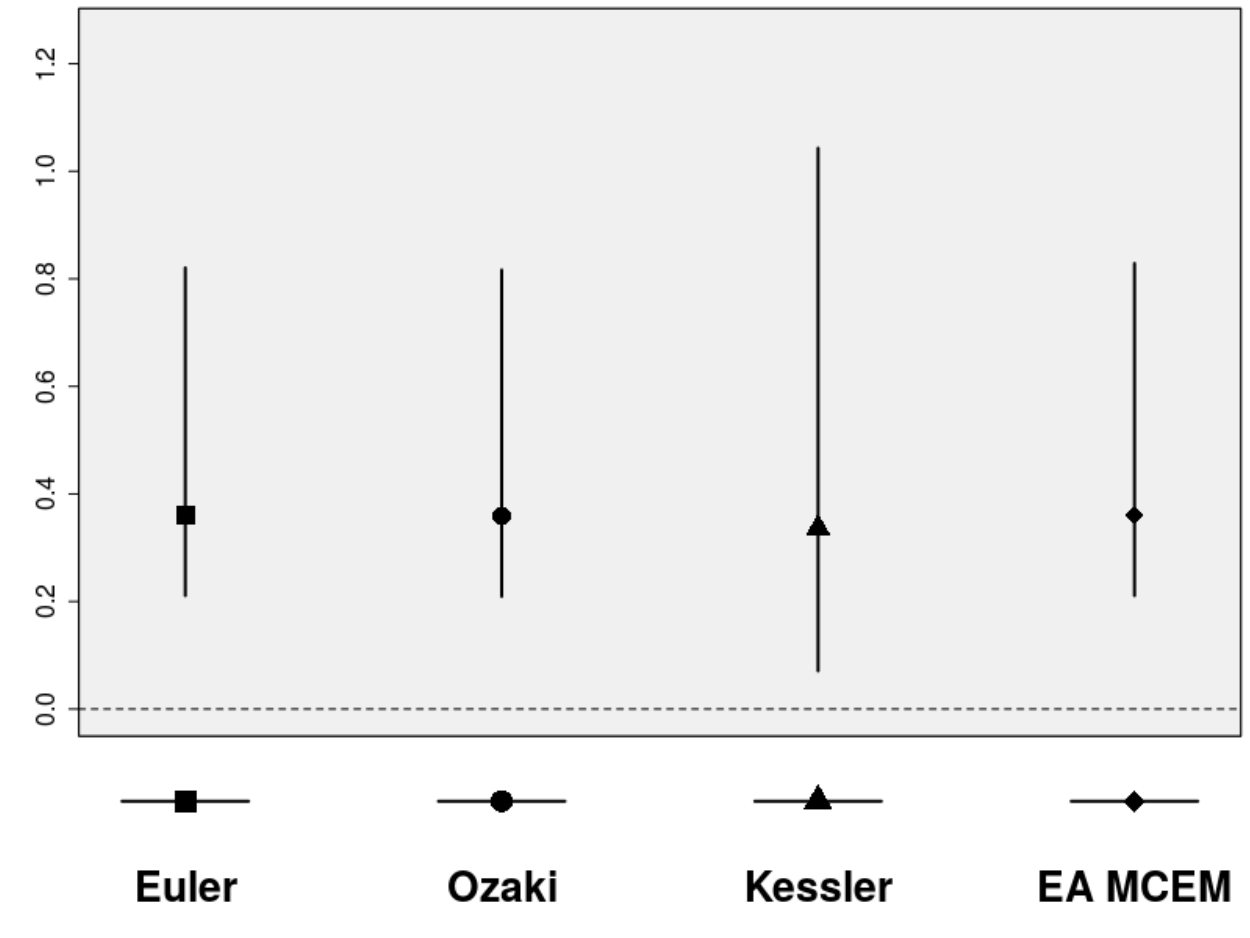}
 \includegraphics[width=0.33\textwidth]{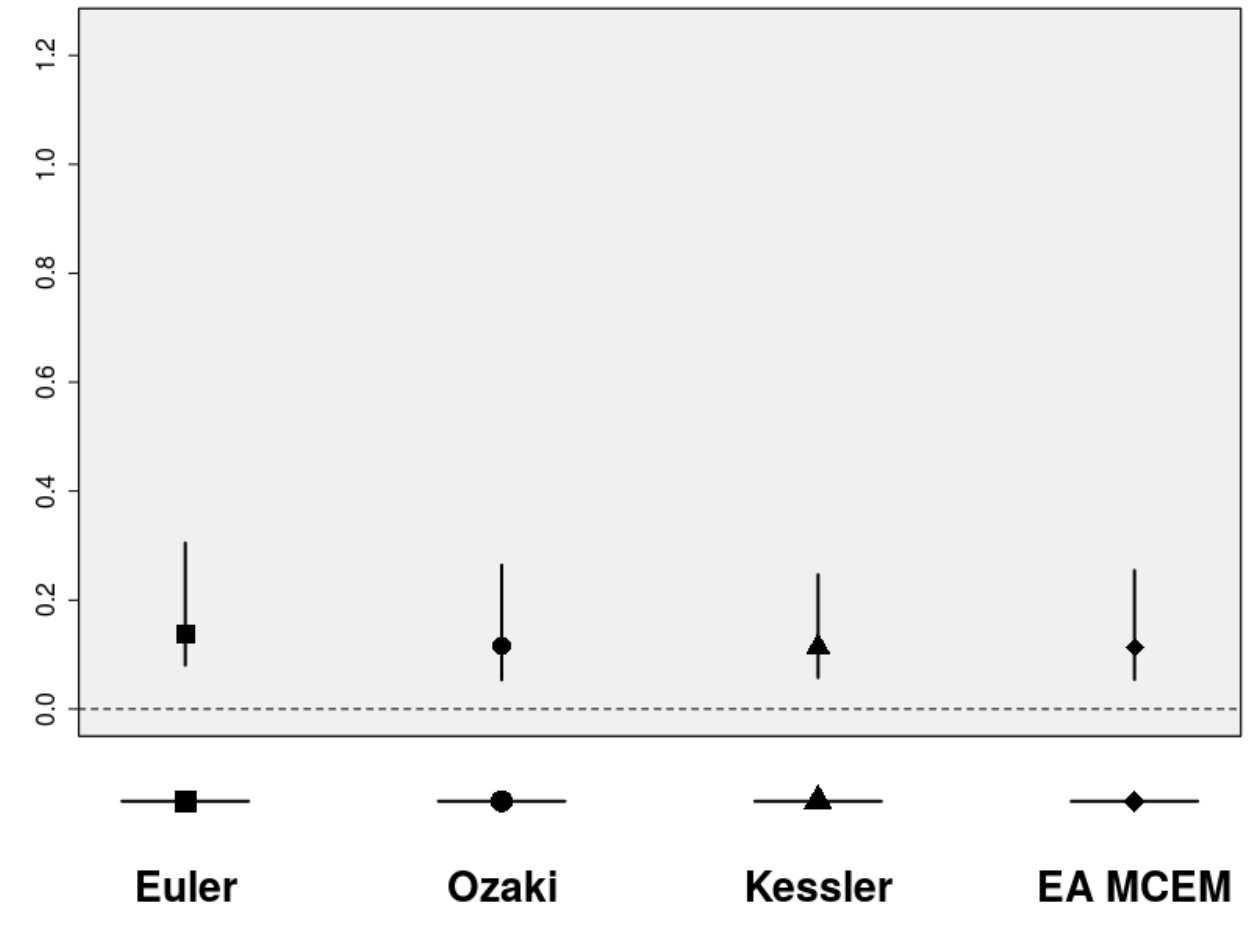}
 \includegraphics[width=0.33\textwidth]{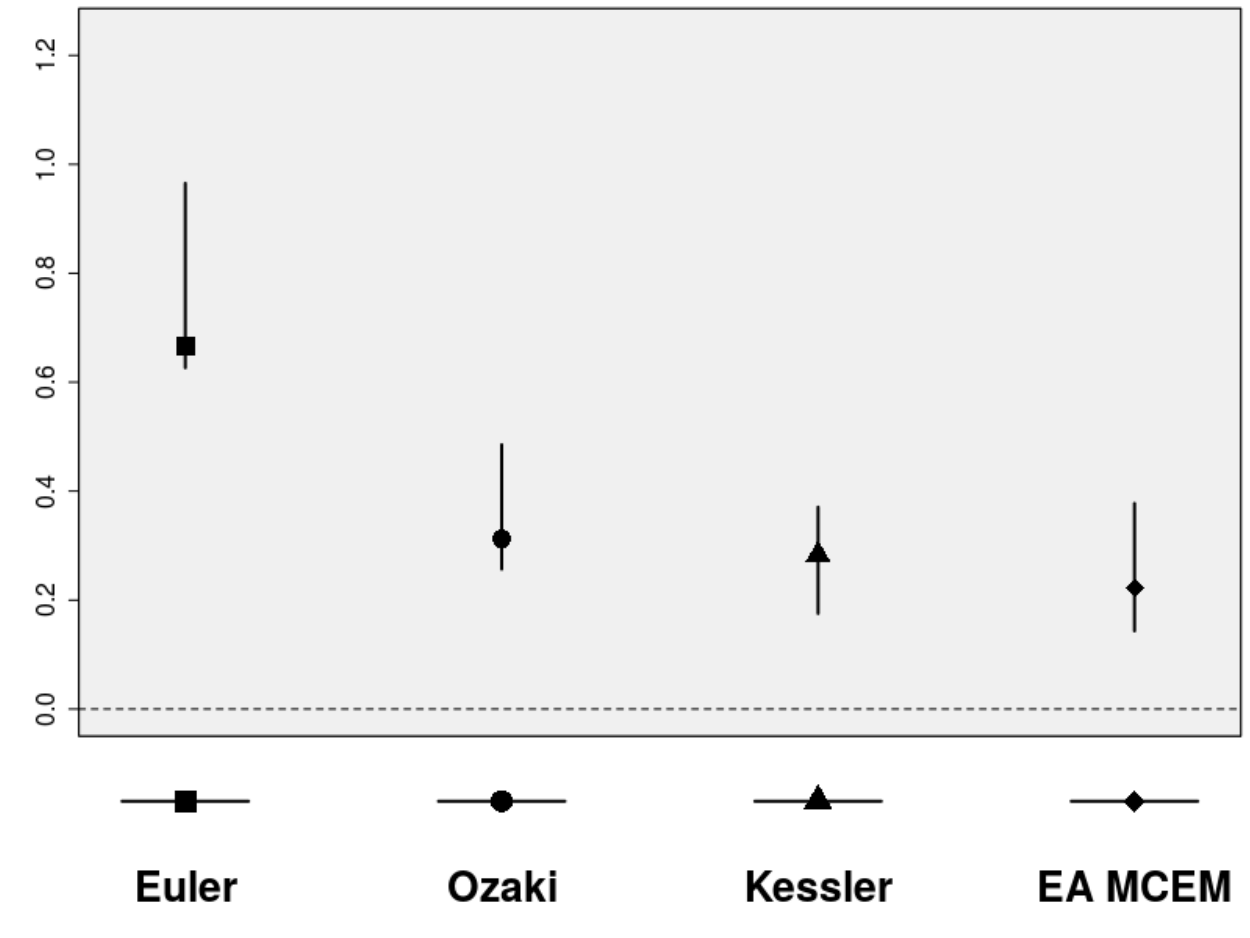}
 \end{tabular}
 \caption{Distribution of the relative integrated error between the estimated map and the true one. Results are presented for each method with $\Delta=0.1$ (left), 1 (center) and 10 (right).}
  \label{fig:sq:error}
 \end{figure}

Figure \ref{fig:sq:error} shows the distribution of the relative integrated error $\frac{\int |P_{\hat{\eta}}(\bo{x}) - P_{\eta^*}(\bo{x})| \rmd \bo{x}}{\int  |P_{\eta^*}(\bo{x}) |\rmd \bo{x}}   $ for each method. 
According to this simulation study,  it seems preferable not to use the Euler method in comparison to any of the three other procedures.  
\subsection{Inferring fishing zones}
\label{subsec:data:application}
\begin{figure}[p]
\centering
\includegraphics[width=0.45\textwidth]{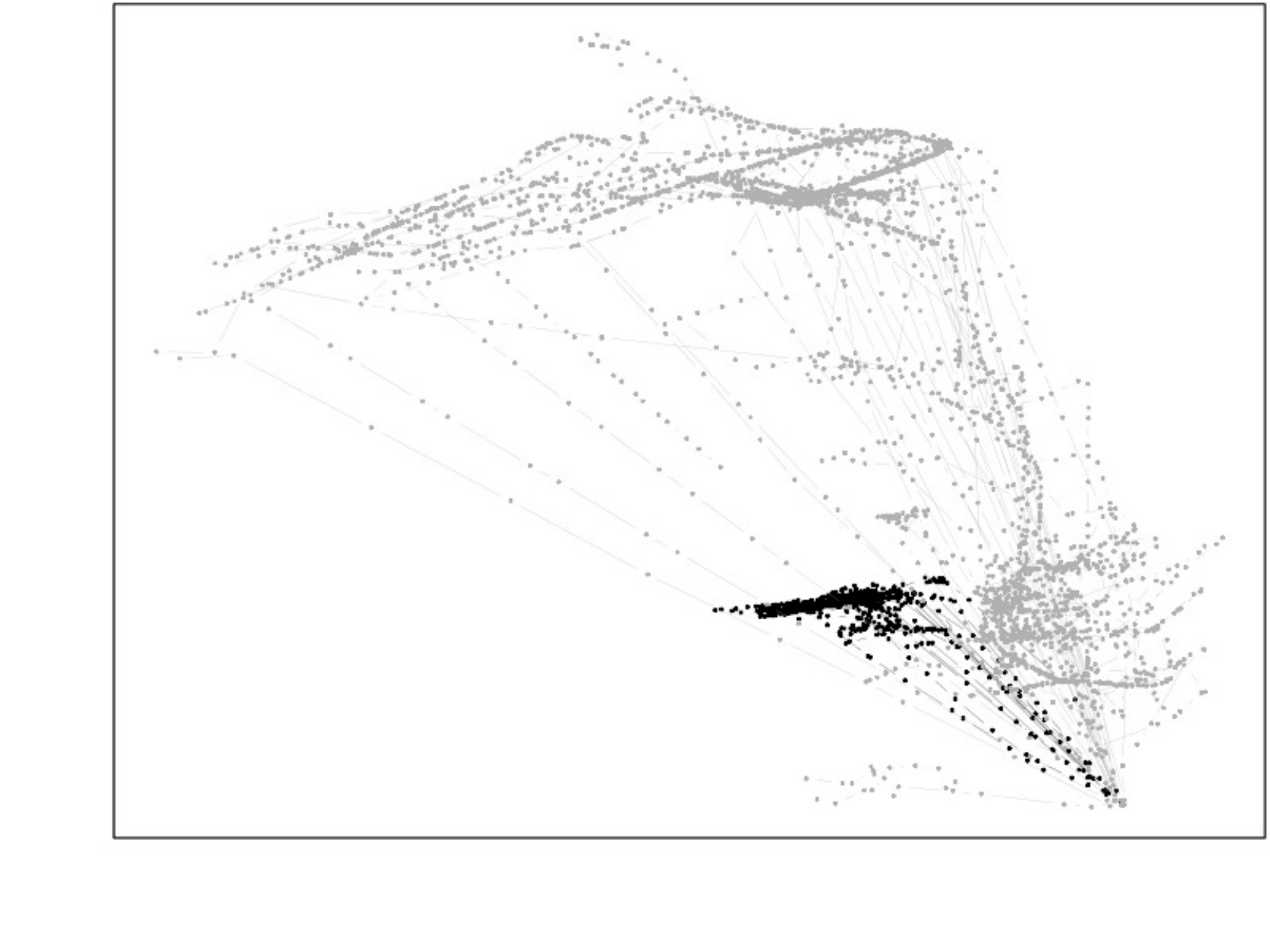}
\includegraphics[width=0.45\textwidth]{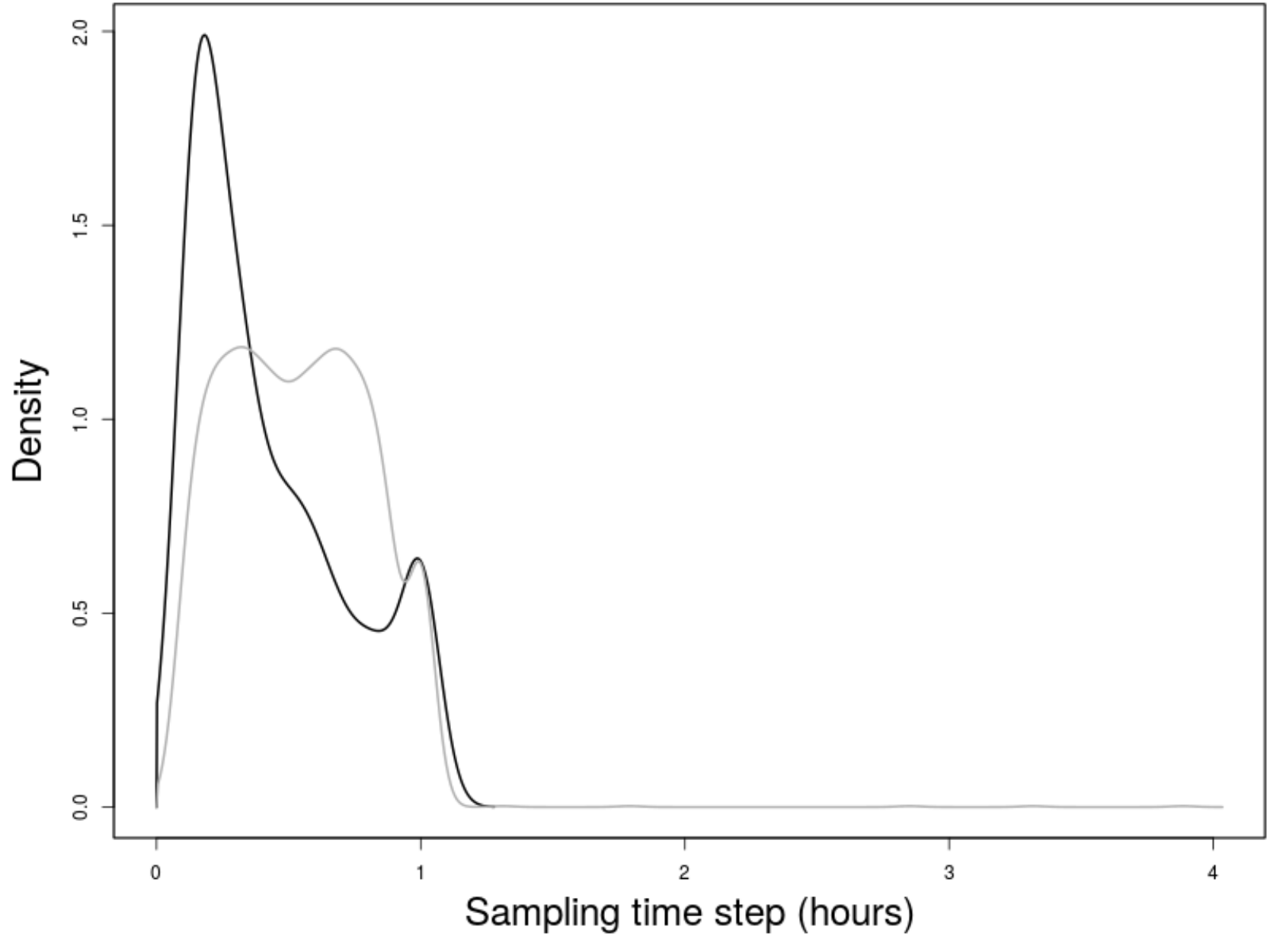}
\caption{\textit{Left}: the two sets of trajectories are presented with a different color, grey or black, depending on the vessel. \textit{Right}: distribution of sampling time step for GPS tags are presented for both datasets. 
 For confidentiality reasons, the actual recorded positions may not be shown and therefore no labels are presented on the axes.}
\label{fig:actual:data}
\end{figure}

Delimiting fishing zones with high potential represents a key step in fisheries management, for instance to set up some marine protected areas. 
In order to illustrate the performance of each method to define such zones using actual GPS data, the model presented in Section \ref{sec:model} is estimated using two sets of trajectories of two French fishing vessels named V1 and V2. 
The two sets of trajectories are represented in Figure \ref{fig:actual:data}.
\begin{itemize}
\item[V1] (Black in Figure \ref{fig:actual:data})  The dataset of vessel 1 is composed of 15 short trajectories with a total of 724 GPS locations and a sampling time step of around 12 minutes, with some irregularities up to 1 hour.
\item[V2] (Grey in Figure \ref{fig:actual:data})  The dataset of vessel 2 is composed of 25 trajectories with a total of 3111 GPS locations and an irregular sampling time step mostly between 15 and 50 minutes, with some irregularities up to 4 hours.  
\end{itemize}
For each dataset, the four methods presented above were used to estimate $\bs{\theta}$. 
The number of components $K$ is not assumed to be known anymore and has to be estimated. We do not focus on this problem which has already been addressed in the case of discretely observed SDE (see \citealp{uchida2005aic}, for instance) and assume $K_1=1$ (respectively $K_2=3$) for V1 (resp. V2). 
In the following, this approximation of the likelihood proposed by \cite{beskos:papaspiliopoulos:roberts::2009} is used to choose which estimation method provides the best estimated map.

\begin{figure}[p]
\centering
\begin{subfigure}{0.49\textwidth}
\caption{Euler}
\includegraphics[width=0.99\textwidth]{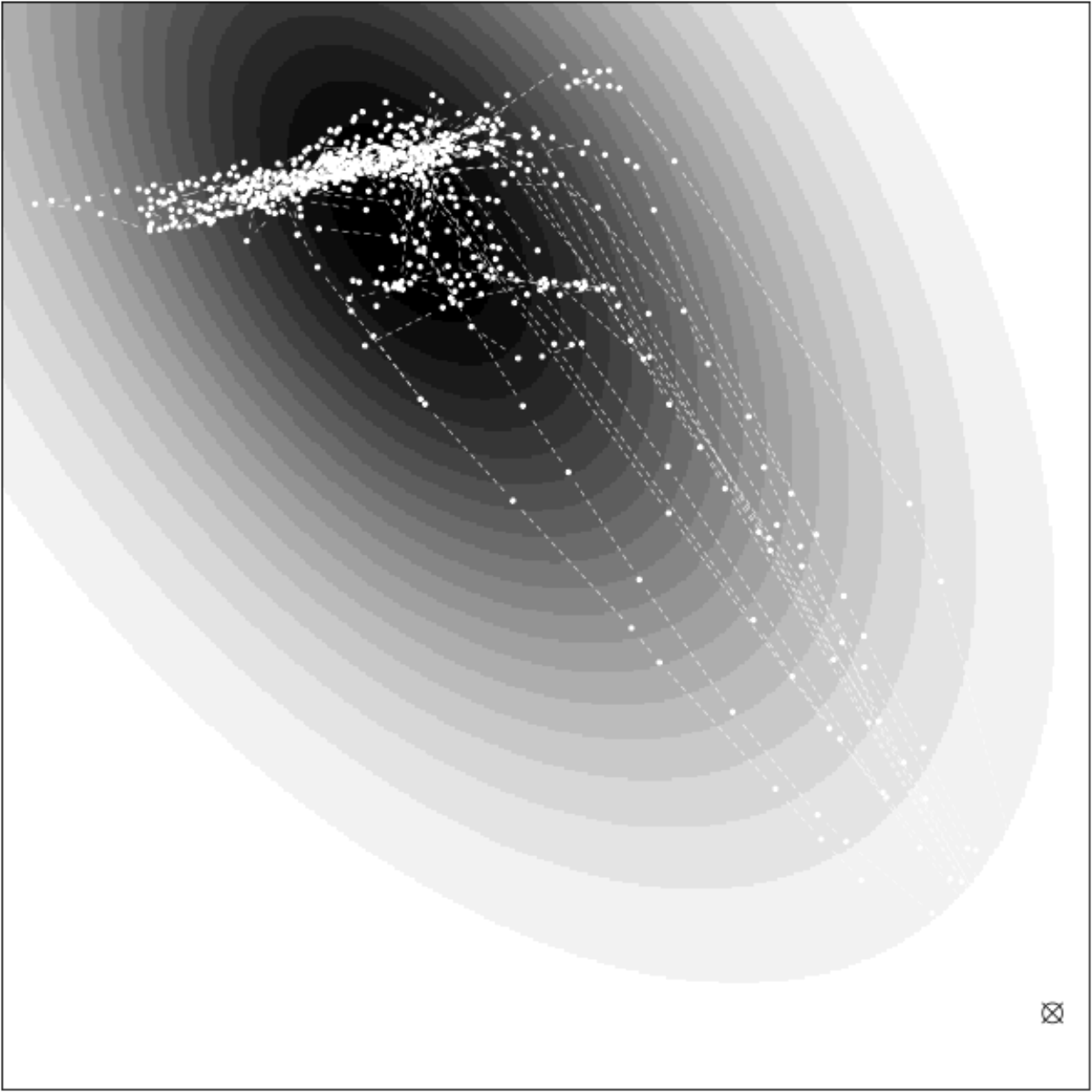}
\end{subfigure}
\begin{subfigure}{0.49\textwidth}
\caption{Kessler}
\includegraphics[width=0.99\textwidth]{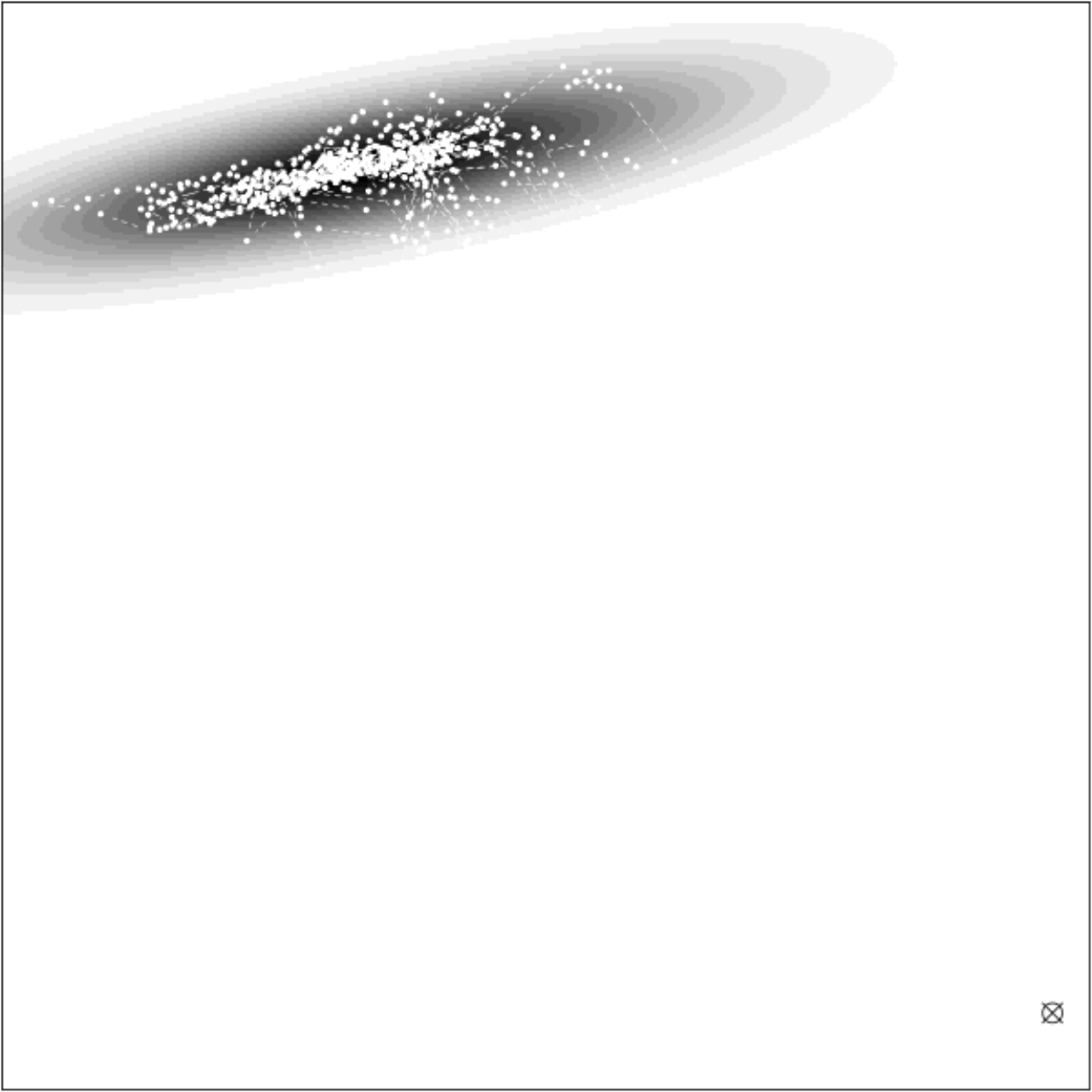}
\end{subfigure}
\begin{subfigure}{0.49\textwidth}
\includegraphics[width=0.99\textwidth]{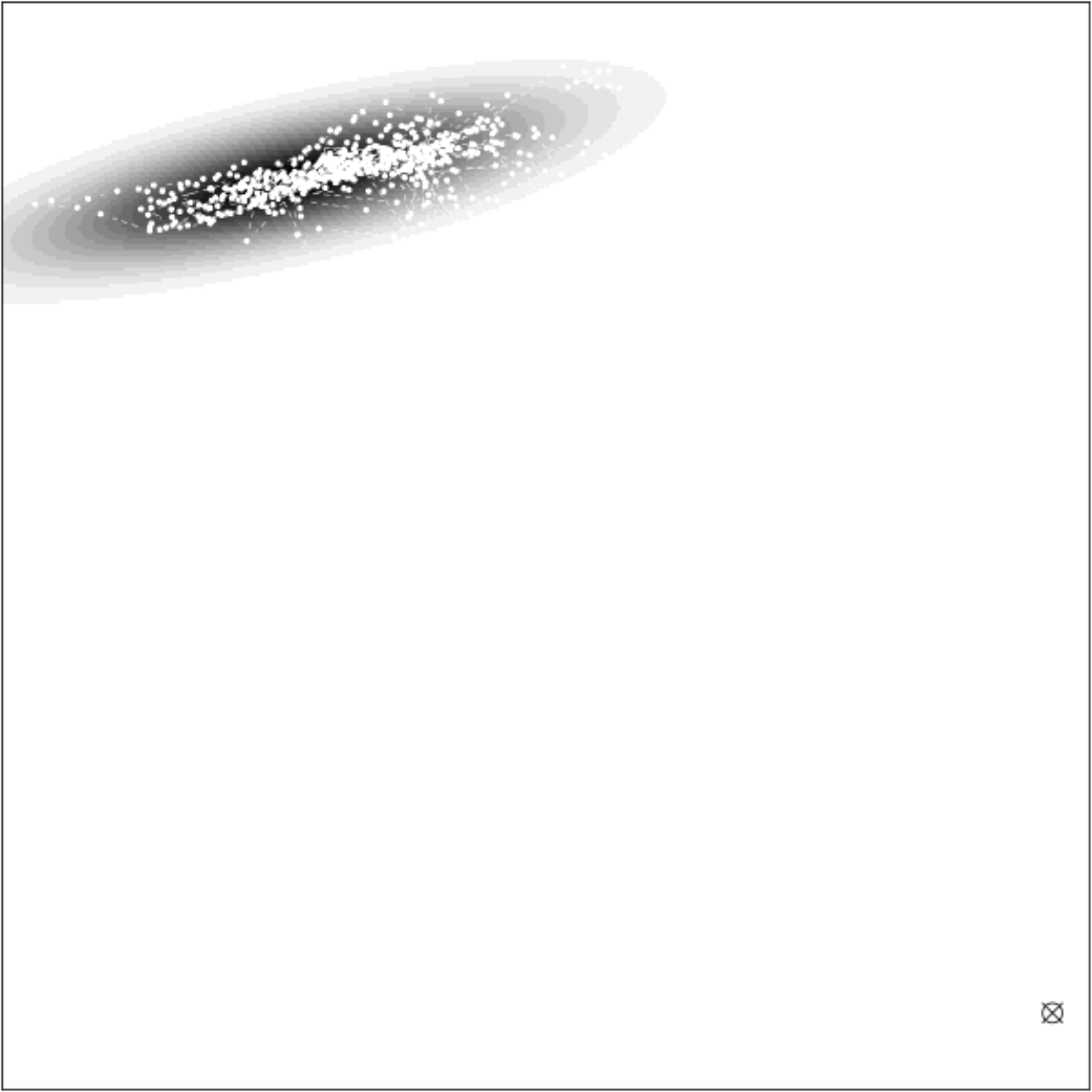}
\caption{Ozaki}
\end{subfigure}
\begin{subfigure}{0.49\textwidth}
\includegraphics[width=0.99\textwidth]{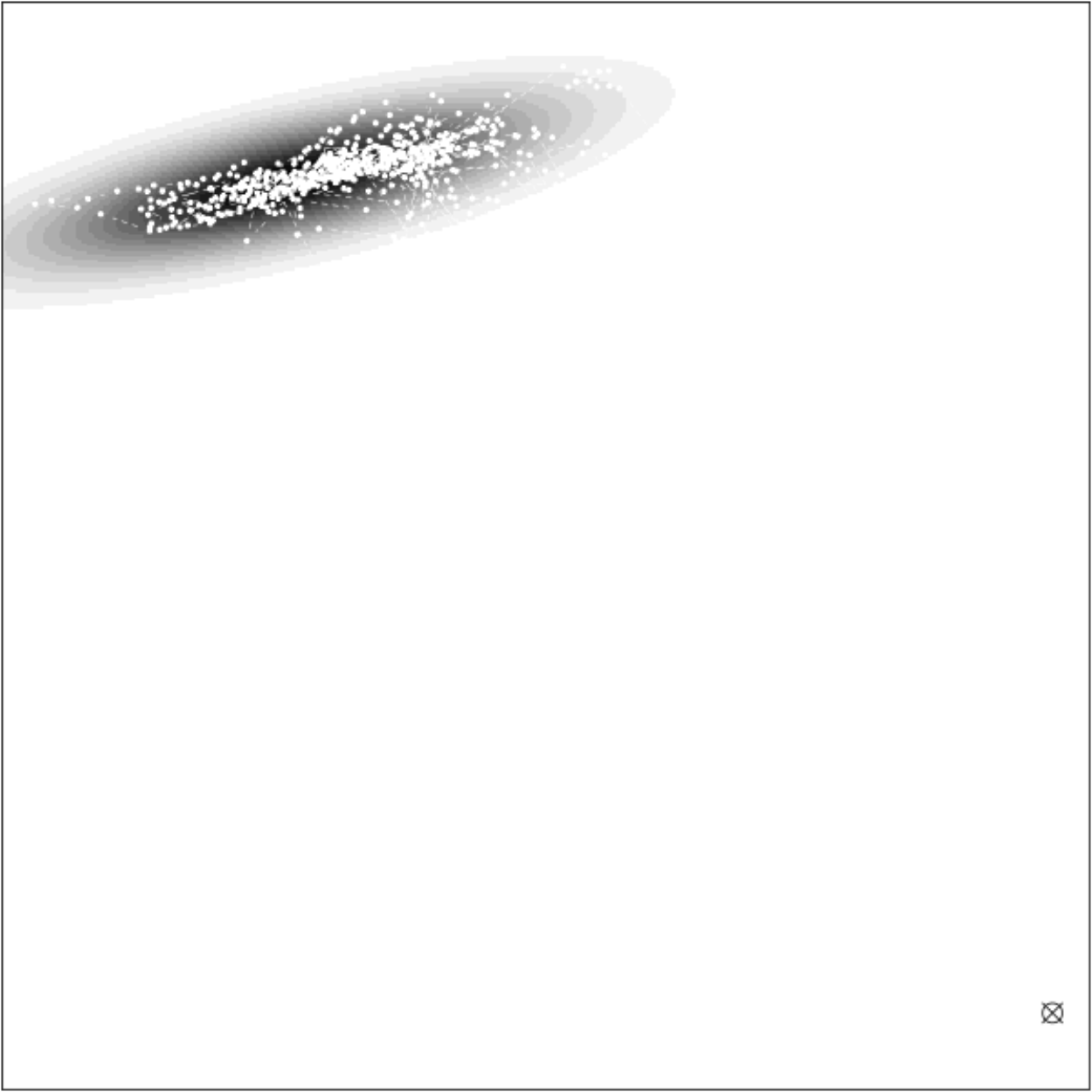}
\caption{EA MCEM}
\end{subfigure}
\caption{Estimated map for the dataset of the vessel 1 (see Figure \ref{fig:actual:data}), using four different estimation methods. 
The cross at the bottom right of each map is the departure harbour. 
The darker a zone is, the more attractive it is for the given vessel. 
Observed points are plotted in white in order to see the superposition between maps and trajectories.}
\label{fig:res:v1}
\end{figure}

\begin{figure}[p]
\centering
\begin{subfigure}{0.49\textwidth}
\caption{Euler}
\includegraphics[width=0.99\textwidth]{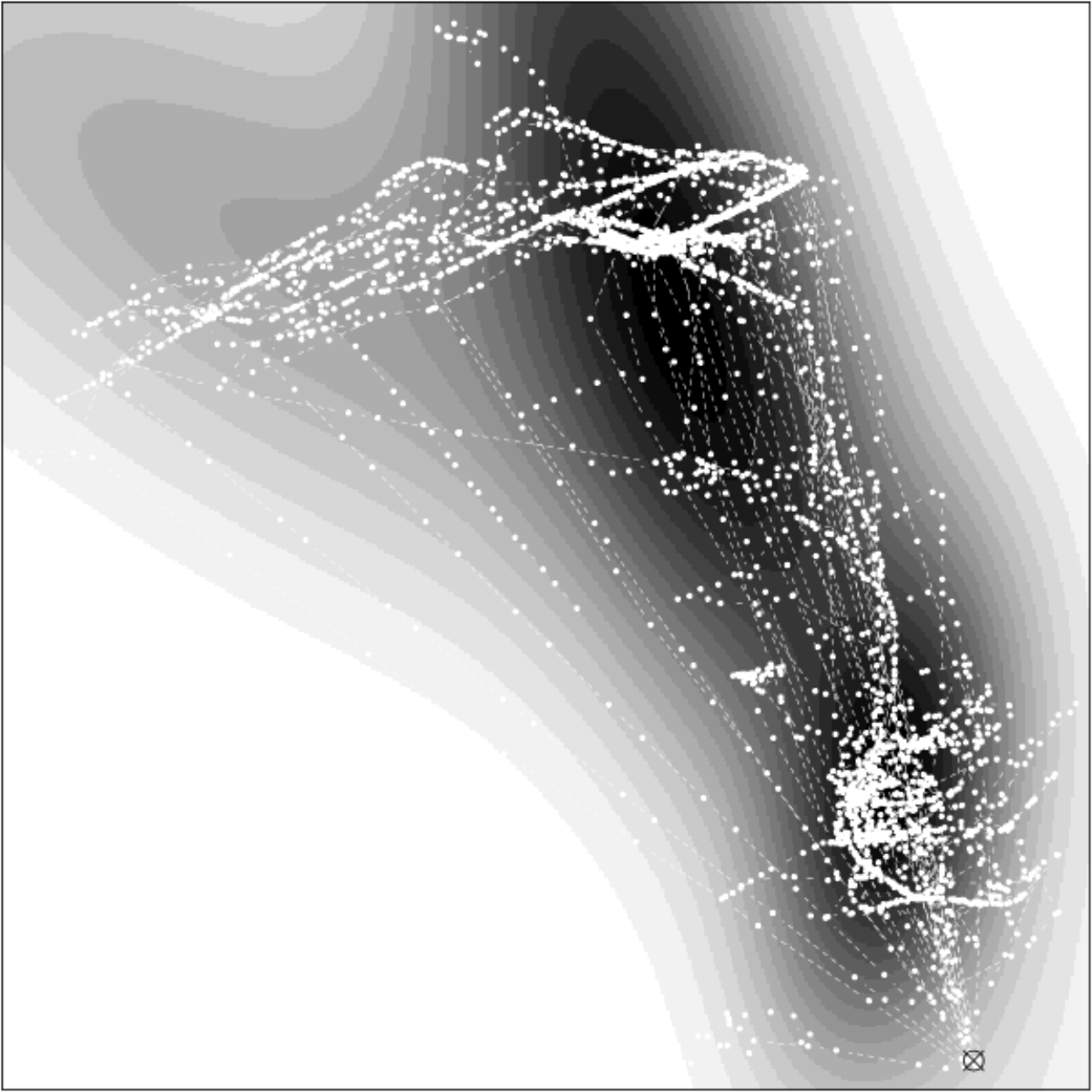}
\end{subfigure}
\begin{subfigure}{0.49\textwidth}
\caption{Kessler}
\includegraphics[width=0.99\textwidth]{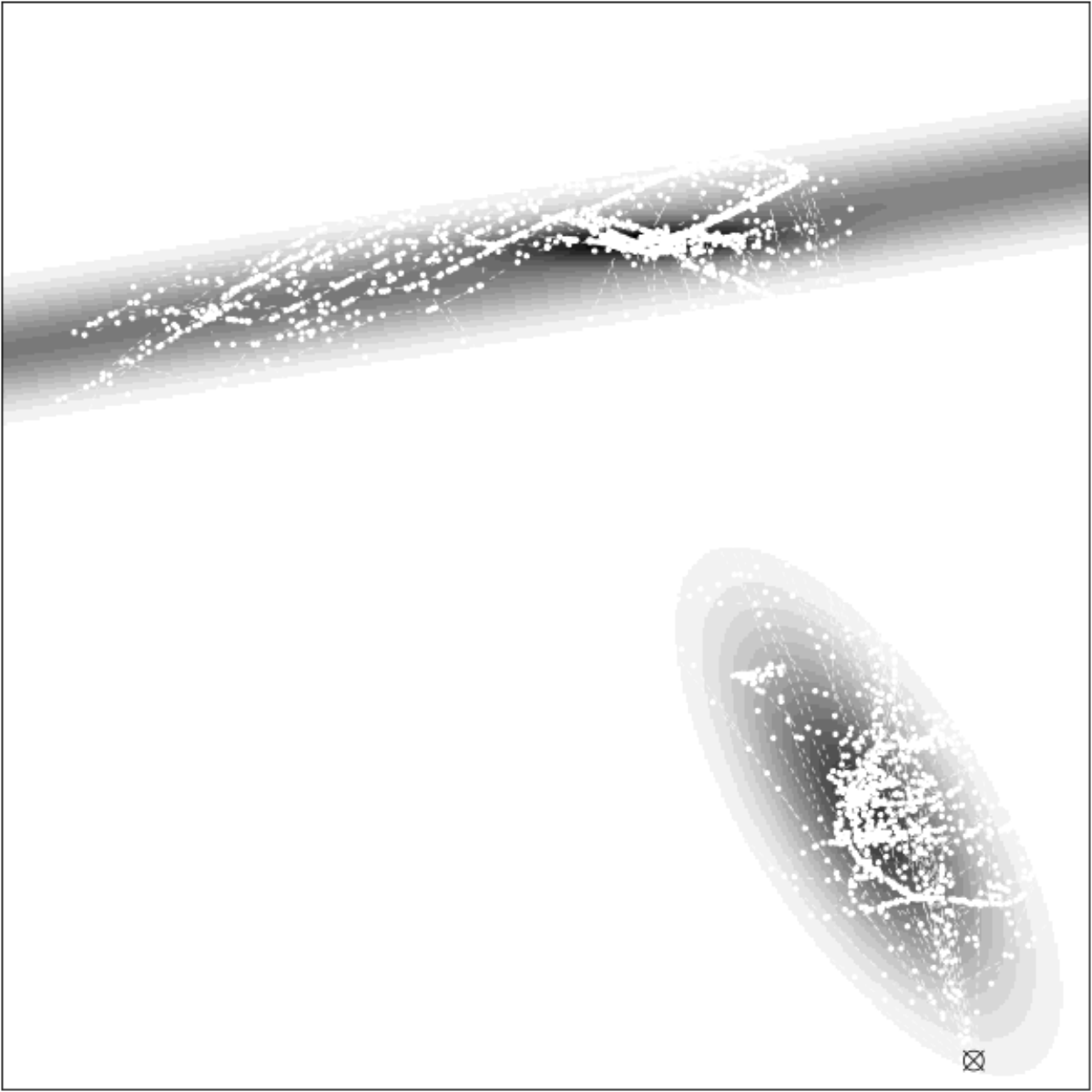}
\end{subfigure}
\begin{subfigure}{0.49\textwidth}
\includegraphics[width=0.99\textwidth]{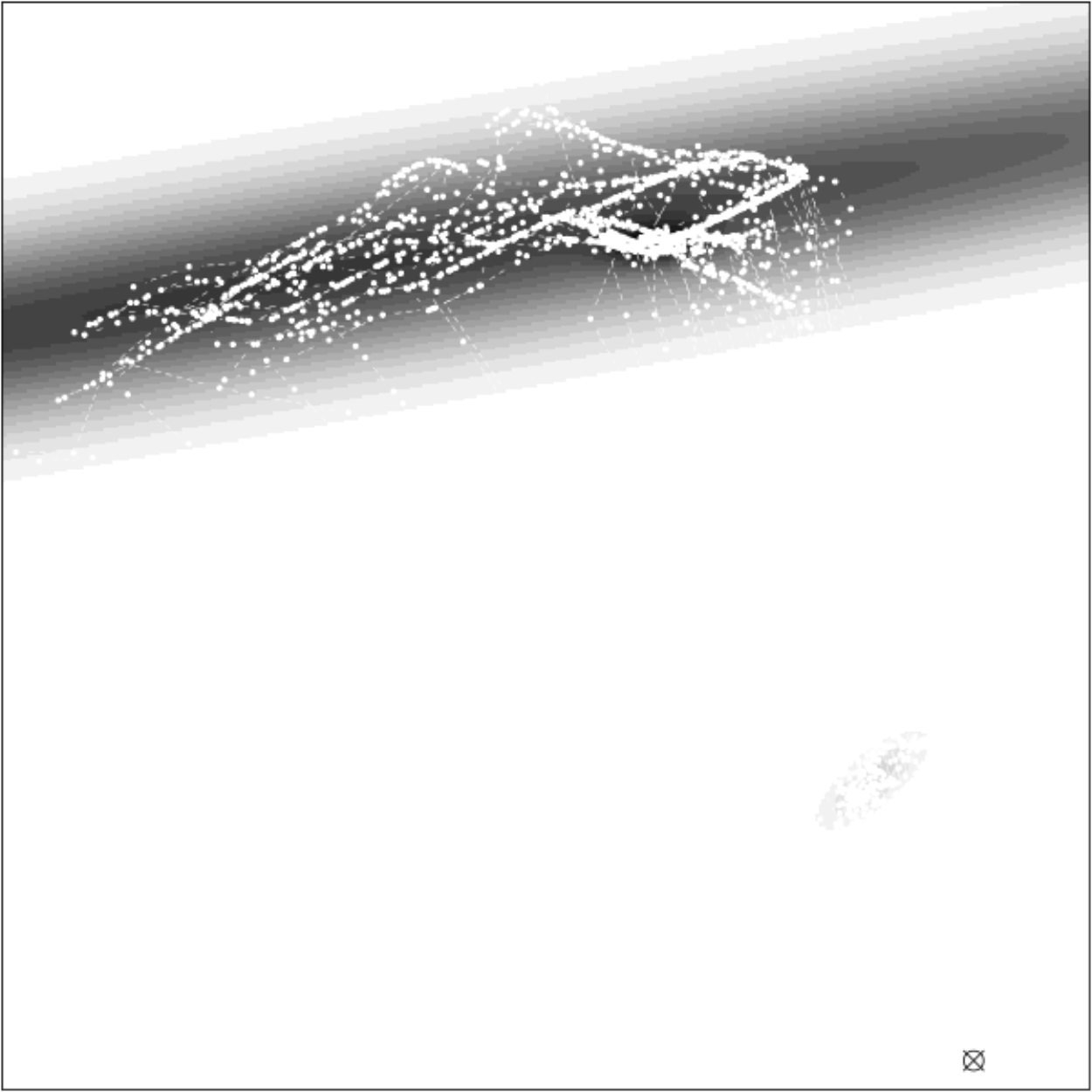}
\caption{Ozaki}
\end{subfigure}
\begin{subfigure}{0.49\textwidth}
\includegraphics[width=0.99\textwidth]{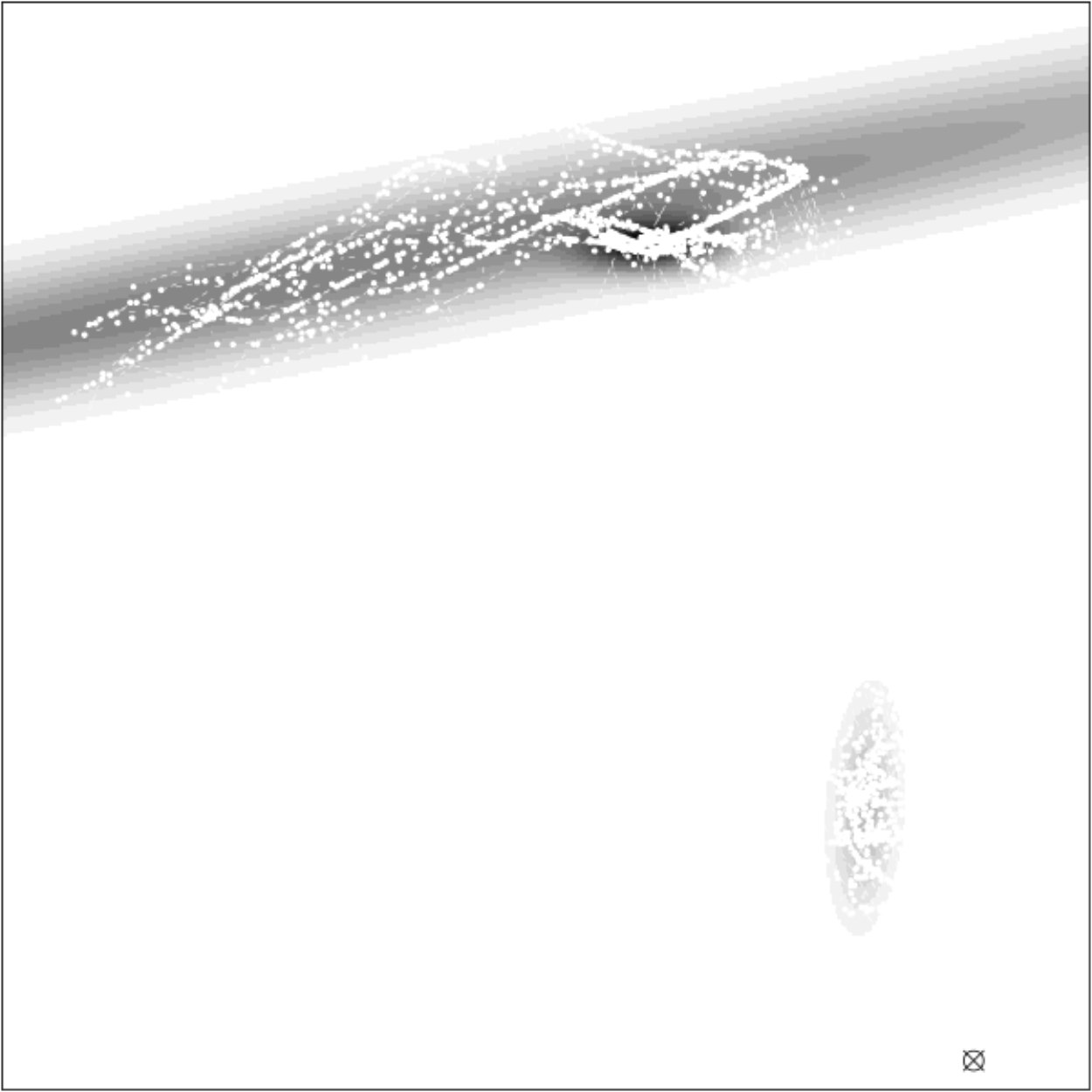}
\caption{EA MCEM}
\end{subfigure}
\caption{Estimated map for the data set of the vessel 2 (see Figure \ref{fig:actual:data},) using four different estimation methods. 
The cross at the bottom right of each map is the departure harbour. 
The darker a zone is, the more attractive it is for the given vessel. 
Observed points are plotted in white in order to see the superposition between maps and trajectories.}
\label{fig:res:v2}
\end{figure}

\paragraph{Results for V1}
Results for the vessel 1 are presented in Figure~\ref{fig:res:v1}. 
The best estimate in the sense of the greatest approximated log-likelihood is given by the  map estimated with the MCEM EA method, followed by Ozaki, Kessler and Euler (see Table \ref{tab:contrast:vals} for the values of each log-likelihood, for each estimated parameter). Again $\hat{\bs{\theta}}_{\mathsf{K}}$ is unstable as the correction term for the variance (given in equation \eqref{eq:kes:var}) leads to a non positive semi-definite matrix for 28\% of the observations. 
However, the estimated map is similar to the one given by Ozaki and EA MCEM methods although less concentrated around the data.  The acceptance rate for the conditional simulation  of the EA MCEM method is 1.4\%.
Euler discretization method estimates a potential map with a much more spread attractive zone than the three other methods. Moreover, the orientation of this zone does not follow what seems to be the main axis of trawling (mainly East-West).  Both Ozaki and MCEM EA methods provide a similar estimated map (in terms of parameters and of likelihood): a Gaussian form wrapping what appears to be the main trawling zone. 
The axis of the resulting Gaussian form are conform to the trawling directions of the vessels.
\paragraph{Results for V2}
Results for the vessel 2 are presented in Figure \ref{fig:res:v2}. As in the case of V1, the best estimate in the sense of the greatest approximated loglikelihood is given by the  map estimated with the EA MCEM method, followed by Ozaki, Kessler and Euler (see Table \ref{tab:contrast:vals} for the values of each log-likelihood, for each estimated parameter). 
Again, $\hat{\bs{\theta}}_{\mathsf{K}}$ is unstable, as the correction term for the variance leads to a non positive semi-definite  matrix for 48\% of observations. The acceptance rate for the conditional simulation  of the EA MCEM method is 6.2\%. As for vessel 1, Euler method produces a smooth estimated map, where attractive zones are connected and wrap almost all observed points. Kessler, Ozaki and EA MCEM methods lead to maps where zones are disconnected: one close to to the harbour and an offshore zone. 
In the case of the Kessler method, the extension and the relative weight of this first zone is larger than for the two other methods. For the three methods, the second offshore zone is a mixture of a general attractive zone with an East-West orientation, and a smaller hot spot that gathers a large amount of observed points. 
Here, a natural extension would be to use a  larger $K$.  The Gaussian model proposed by \cite{preisler:et:al:2013} would have imposed a potential function whose orientation is given by the x and y axis, with circular contour lines as a consequence of the independence between directions. 
The generalization proposed in this paper models more complex attractive zones. 
The contour lines of the potential might be used to define more realistic high potential fishing zones.

\begin{alternativetable}
\centering
\begin{tabular}{|c|cccc|cccc|} 
\multicolumn{1}{c}{}&\multicolumn{4}{c}{V1}& \multicolumn{4}{|c|}{V2}\\
\cline{2-9}
\multicolumn{1}{c|}{Estimate  $\backslash$ Contrast} & $\ell_\mathsf{E}$ & $\ell_\mathsf{O}$ & $\ell_\mathsf{K}$ & $\ell_\mathsf{EA}$& $\ell_\mathsf{E}$ & $\ell_\mathsf{O}$ & $\ell_\mathsf{K}$ & $\ell_\mathsf{EA}$ \\ 
  \hline
$\bs{\theta}_\mathsf{E}$   & -1.32 & -1.29 & -0.37 & -1.29 & -2.17 & -2.16 & -1.07 & -2.16 \\ 
 $\bs{\theta}_\mathsf{O}$  & -1.50 & -0.91 &  0.39 & -0.93 & -2.21 & -2.09 & -0.90 & -2.09 \\ 
  $\bs{\theta}_\mathsf{K}$ & -1.50 & -0.96 &  0.70 & -0.99 & -2.26 & -2.18 & -0.67 & -2.10 \\ 
  $\bs{\theta}_\mathsf{EA}$  & -1.49 & -0.92 &  0.06 & -0.92 & -2.21 & -2.11 & -0.98 & -2.07 \\ 
   \hline
\end{tabular}
\caption{\label{tab:contrast:vals}
Value of the contrast function to maximize for each estimated parameter, for each data set. Each entry is the value of the contrast function (given in column) evaluated for an estimation of $\bs{\theta}$ (given in line). Highest value for each column is on the diagonal. Values are normalized contrasts (divided by the number of segments of trajectories).}
\end{alternativetable}
\section{Conclusions}
This paper proposes a new model based on a stochastic differential equation to describe movement in ecology. 
The proposed model belongs to the popular class of continuous time model described by SDE, whose drift is supposed to be the gradient of a potential function. 
This potential function is assumed to reflect the attractiveness of the environment in which the moving individual evolves.
 We propose a model where this potential is a mixture of smooth Gaussian shaped functions, giving a flexible framework to describe multimodal potential surfaces. As a counterpart, this model leads to a difficult inference task, as some key quantities to express the likelihood are unknown. This lead us explore different possible approaches for parameters estimation.

Several extensions of the proposed potential based model could be considered to capture sharply the individuals dynamics. 
For instance, in this paper, the first part of the fishing vessels trajectories, which consists in steaming toward a fishing zone, seems to be more deterministic than the second part, when the vessel is in the fishing zone. 
We could use a Markov switching SDE where the state might follow several dynamics depending on a regime indicator given by a finite state space hidden Markov chain.  
In the case presented here, this deterministic phase contains a small proportion of the trajectory, and therefore does not affect the inference result, but a better model for this first phase might be preferable.

We also make the assumptions that observed trajectories are independent. 
In many cases, this is likely to be false. Therefore, the $G$ observed trajectories could also be modeled as the solution to a $2G$ dimensional SDE with interactions in the drift function and the diffusion matrix. 
As the individuals are not observed at the same time, a partially observed SDE can be defined which implies additional challenges to estimate all parameters using algorithms designed for hidden Markov models. 

Finally, an undesirable property of the potential introduced in the paper is that, far from the attractive zones, the process behaves like a Brownian motion. 
This can be easily overcome for instance by adding another component in the mixture proportional to $-\|\bo{x}-\mathsf{G}\|^{1/2}$, where $\mathsf{G}$ is a fixed attractive point. 
This ensures that even if an individual is far from the attractive zones, he is likely to revert around $\mathsf{G}$.  
We used this potential and obtained very similar numerical results. 
We decided not to add this component in the model as it seemed artificial to ensure stationarity for theoretical reasons.
\vspace{\baselineskip}

From the inference point of view, this work explores the performance of four inference methods to estimate the parameters in such potential based movement models. 
The first three methods are pseudo likelihood approaches, including the widely used Euler-Maruyama procedure, with similar implementations and requiring little assumptions on the model.
The fourth method, based on a  Monte Carlo EM using  the Exact Algorithm, requires additional mathematical assumptions on the model and is more technically and computationally intensive, but has the nice property of not using Gaussian approximations, inducing a structural bias.

Numerical experiments with simulated data illustrate that the Euler discretization method is not robust to low sampling rates, which is a common situation in movement ecology, and that the three other approximate maximum likelihood methods perform better, suggesting that the Euler approach might not be well suited to movement ecology.
 
These four methods were also compared to estimate potential maps (from a specific model) using actual movement data of two fishing vessels. 
It is clear here that the different methods can lead to different maps with various interpretations. 
Specifically, in the examples presented here, we would like to point out that the Euler-Maruyama results lead to a totally different map, implying different interpretations, than the three other methods. 
Therefore, users should be careful when using this method, and ensures that their conclusions are not biased by the choice of the Euler discretization. 
Moreover, using the criterion of the unbiased approximated likelihood of \cite{beskos:papaspiliopoulos:roberts::2009}, the Euler method shows the worst results. 

Among the alternative likelihood approaches the Kessler approximate method was more unstable than the Ozaki discretization as the proposed approximation of the  covariance matrix can be non positive semi-definite for some observed points. 
The Kessler method is only used in this work with an expansion of conditional moments up to order two and higher order expansions could lead to more robust estimates.
The Ozaki method showed results similar to the  unbiased EA MCEM approach, which were the best according to the approximated likelihood criterion.

Overall, the EA MCEM based method seems to be the most appealing as it does not introduce any bias, and, in our study, was much less sensitive to the starting point than the pseudo likelihood methods. 
However, the difficulty of its implementation and the computation cost might reduce its interest. 
Moreover, its application is restricted to a specific kind of potential functions. 

Based on simulations and the real dataset presented in this paper, the Ozaki method seems to provide close estimate to the one given by the EA MCEM method, and is as easy as the Euler approach to implement. 
In addition, it requires only weak assumptions (the drift must have an invertible Jacobian matrix at observed points) and, unlike the Kessler approach, this method is numerically stable. 
Therefore, we suggest to users in movement ecology the use of the Ozaki method to estimate parameters from movement data.

\paragraph{Acknowledgements:} The authors would like to thank the associate editor and the anonymous referees for their meaningful  comments which led to a major reorientation of this work. 
Meeting expenses between the three authors have been supported by the GdR EcoStat. 
The real data application would not have been possible without all contributors to the RECOPESCA project. 
The authors express their warm thanks to the leaders of this project Patrick Berthou and \'Emilie Leblond and to all the voluntary fishermen. 
The authors have also benefited from insightful discussions with Valentine Genon-Catalot and Catherine Laredo.
\bibliographystyle{rss}
\bibliography{biblio}

\begin{thebibliography}{32}
\expandafter\ifx\csname natexlab\endcsname\relax\def\natexlab#1{#1}\fi
\expandafter\ifx\csname url\endcsname\relax
  \def\url#1{\texttt{#1}}\fi
\expandafter\ifx\csname urlprefix\endcsname\relax\def\urlprefix{URL: }\fi

\bibitem[{Ait-Sahalia(1999)}]{ait-sahalia:1999}
Ait-Sahalia, Y. (1999) Transition densities for interest rate and other
  nonlinear diffusions.
\newblock \textit{Journal of {Finance}}, \textbf{54}, 1361--1395.

\bibitem[{Ait-Sahalia(2002)}]{ait-sahalia:2002}
--- (2002) Maximum likelihood estimation of discretely sampled diffusions: a
  closed-form approximation approach.
\newblock \textit{Econometrica}, \textbf{70}, 223--262.

\bibitem[{Ait-Sahalia(2008)}]{ait-sahalia:2008}
--- (2008) Closed-form likelihood expansions for multivariate diffu- sions.
\newblock \textit{The {A}nnals of {S}tatistics}, \textbf{36}, 906--937.

\bibitem[{Beskos et~al.(2006{\natexlab{a}})Beskos, Papaspiliopoulos and
  Roberts}]{beskos:papaspiliopoulos:roberts:2006}
Beskos, A., Papaspiliopoulos, O. and Roberts, G. (2006{\natexlab{a}})
  {Retrospective exact simulation of diffusion sample paths with applications}.
\newblock \textit{Bernoulli}, \textbf{12}, 1077--1098.

\bibitem[{Beskos et~al.({2009})Beskos, Papaspiliopoulos and
  Roberts}]{beskos:papaspiliopoulos:roberts::2009}
--- ({2009}) Monte {C}arlo maximum likelihood estimation for discretely
  observed diffusions processes.
\newblock \textit{{Annals of {S}tatistics}}, \textbf{{37}}, {223--245}.

\bibitem[{Beskos et~al.(2006{\natexlab{b}})Beskos, Papaspiliopoulos, Roberts
  and Fearnhead}]{beskos:papaspiliopoulos:roberts:fearnhead:2006}
Beskos, A., Papaspiliopoulos, O., Roberts, G. and Fearnhead, P.
  (2006{\natexlab{b}}) {Exact and computationally efficient likelihood-based
  estimation for discretely observed diffusion processes}.
\newblock \textit{Journal of Royal Statistical Society, Series B: Statistical
  Methodology}, \textbf{68}, 333--382.

\bibitem[{Beskos and Roberts(2005)}]{beskos:roberts:2005}
Beskos, A. and Roberts, G. (2005) {Exact simulation of diffusions}.
\newblock \textit{Annals of {Applied} {P}robability}, \textbf{15}, 2422--2444.

\bibitem[{Blackwell(1997)}]{blackwell:1997}
Blackwell, P. (1997) Random diffusion models for animal movement.
\newblock \textit{Ecological {M}odelling}, \textbf{100}, 87 -- 102.
\newblock
  \urlprefix\url{http://www.sciencedirect.com/science/article/pii/S0304380097001531}.

\bibitem[{Blackwell et~al.(2015)Blackwell, Niu, Lambert and
  LaPoint}]{blackwell:et:al:2015}
Blackwell, P., Niu, M., Lambert, M. and LaPoint, S. (2015) Exact {B}ayesian
  inference for animal movement in continuous time.
\newblock \textit{Methods in {E}cology and {E}volution}.

\bibitem[{Brillinger(2010)}]{brillinger:2010}
Brillinger, D. (2010) \textit{Handbook of Spatial Statistics}, chap.~26.
\newblock Chapman and Hall/CRC Handbooks of Modern Statistical Methods. CRC
  Press.

\bibitem[{Brillinger et~al.(2002)Brillinger, Haiganoush, Ager, Kie and
  Stewart}]{brillinger:et:al:2002}
Brillinger, D., Haiganoush, K., Ager, A., Kie, J. and Stewart, B. (2002)
  {Employing stochastic differential equations to model wildlife motion}.
\newblock \textit{Bulletin of the {B}razilian {M}athematical {S}ociety},
  \textbf{33}, 385--408.

\bibitem[{Brillinger et~al.(2001{\natexlab{a}})Brillinger, Preisler, Ager and
  Kie}]{brillinger:et:al:2001a}
Brillinger, D., Preisler, H., Ager, A. and Kie, J. (2001{\natexlab{a}}) The use
  of potential functions in modeling animal movement.
\newblock \textit{Data analysis from statistical foundations}, 369--386.

\bibitem[{Brillinger et~al.(2001{\natexlab{b}})Brillinger, Preisler, Ager, Kie
  and Stewart}]{brillinger:et:al:2001b}
Brillinger, D., Preisler, H., Ager, A., Kie, J. and Stewart, B.
  (2001{\natexlab{b}}) Modeling movements of free-ranging animals.
\newblock \textit{Univ. Calif. Berkeley Statistics Technical Report},
  \textbf{610}.

\bibitem[{Brillinger et~al.(2011)Brillinger, Preisler and
  Wisdom}]{brillinger:et:al:2011}
Brillinger, D., Preisler, H. and Wisdom, M. (2011) Modeling particles moving in
  a potential field with pairwise interactions and an application.
\newblock \textit{Brazilian Journal of Probability and Statistics},
  \textbf{25}, 421--436.

\bibitem[{Chang(2011)}]{Chang2011}
Chang, S.-K. (2011) {Application of a vessel monitoring system to advance
  sustainable fisheries management---Benefits received in Taiwan}.
\newblock \textit{Marine Policy}, \textbf{35}, 116--121.
\newblock
  \urlprefix\url{http://www.sciencedirect.com/science/article/pii/S0308597X10001491}.

\bibitem[{Chavez(2006)}]{Chavez+2006}
Chavez, A. S. G. E.~M. (2006) {Landscape Use and Movements of Wolves in
  Relation to Livestock in a Wildland--Agriculture Matrix}.
\newblock \textit{Journal of Wildlife Management}, \textbf{70}, 1079--1086.

\bibitem[{Dempster et~al.(1977)Dempster, Laird and
  Rubin}]{dempster:maximum:1977}
Dempster, A.~P., Laird, N.~M. and Rubin, D.~B. (1977) {Maximum likelihood from
  incomplete data via the {EM} algorithm}.
\newblock \textit{J. Roy. Statist. Soc. B}, \textbf{39}, 1--38 (with
  discussion).

\bibitem[{Florens-zmirou(1989)}]{florens1989approximate}
Florens-zmirou, D. (1989) Approximate discrete-time schemes for statistics of
  diffusion processes.
\newblock \textit{Statistics}, \textbf{20}, 547--557.
\newblock \urlprefix\url{http://dx.doi.org/10.1080/02331888908802205}.

\bibitem[{Hansen and Ostermeier(2001)}]{hansen_evolutionary_2001}
Hansen, N. and Ostermeier, A. (2001) {Completely derandomized self-adaptation
  in evolution strategies}.
\newblock \textit{Evolutionary {C}omputation}, \textbf{9}, 159--195.

\bibitem[{Harris and Blackwell(2013)}]{harris:blackwell:2013}
Harris, K.~J. and Blackwell, P.~G. (2013) Flexible continuous-time modeling for
  heterogeneous animal movement.
\newblock \textit{Ecological Modelling}, \textbf{255}, 29--37.

\bibitem[{Iacus(2009)}]{iacus2009simulation}
Iacus, S.~M. (2009) \textit{Simulation and inference for stochastic
  differential equations: with R examples}, vol.~1.
\newblock Springer Science \& Business Media.

\bibitem[{Kessler(1997)}]{kessler:1997}
Kessler, M. (1997) Estimation of an ergodic diffusion from discrete
  observations.
\newblock \textit{Scandinavian Journal of Statistics}, \textbf{24}, 211--229.

\bibitem[{Kessler et~al.(2012)Kessler, Lindner and
  Sorensen}]{kessler:lindner:sorensen:2012}
Kessler, M., Lindner, A. and Sorensen, M. (2012) \textit{Statistical methods
  for stochastic differential equations}.
\newblock CRC Press.

\bibitem[{Li(2013)}]{li:2013}
Li, C. (2013) Maximum-likelihood estimation for diffusion processes via
  closed-form density expansions.
\newblock \textit{The {A}nnals of {S}tatistics}, \textbf{41}, 1350--1380.

\bibitem[{Ozaki(1992)}]{ozaki:1992}
Ozaki, T. (1992) A bridge between nonlinear time series models and nonlinear
  stochastic dynamical systems: a local linearization approach.
\newblock \textit{Statistica {S}inica}, 113--135.

\bibitem[{Preisler et~al.(2004)Preisler, Ager, Johnson and
  Kie}]{preisler:et:al:2004}
Preisler, H., Ager, A., Johnson, B. and Kie, J. (2004) Modeling animal
  movements using stochastic differential equations.
\newblock \textit{Environmetrics}, \textbf{15}, 643--657.
\newblock \urlprefix\url{http://dx.doi.org/10.1002/env.636}.

\bibitem[{Preisler et~al.(2013)Preisler, Ager and Wisdom}]{preisler:et:al:2013}
Preisler, H., Ager, A. and Wisdom, M. (2013) Analyzing animal movement patterns
  using potential functions.
\newblock \textit{Ecosphere}, \textbf{4}.

\bibitem[{Shoji and Ozaki(1998{\natexlab{a}})}]{shoji:1998}
Shoji, I. and Ozaki, T. (1998{\natexlab{a}}) Estimation for nonlinear
  stochastic differential equations by a local linearization method 1.
\newblock \textit{Stochastic {A}nalysis and {A}pplications}, \textbf{16},
  733--752.

\bibitem[{Shoji and Ozaki(1998{\natexlab{b}})}]{shojiOzaki1998}
--- (1998{\natexlab{b}}) A statistical method of estimation and simulation for
  systems of stochastic differential equations.
\newblock \textit{Biometrika}, \textbf{85}, 240--243.

\bibitem[{Skellam(1951)}]{skellam_random_1951}
Skellam, J. (1951) Random dispersal in theoretical populations.
\newblock \textit{Biometrika}, 196--218.

\bibitem[{Uchida and Yoshida(2005)}]{uchida2005aic}
Uchida, M. and Yoshida, N. (2005) {AIC} for ergodic diffusion processes from
  discrete observations.
\newblock \textit{preprint MHF}, \textbf{12}.

\bibitem[{Uchida and Yoshida(2012)}]{uchida:yoshida:2012}
--- (2012) Adaptive estimation of an ergodic diffusion process based on sampled
  data.
\newblock \textit{Stochastic Processes and their Applications}, \textbf{122},
  2885 -- 2924.

\end{thebibliography}
\end{document}